\documentclass[12pt]{article}
\usepackage[utf8]{inputenc}
\usepackage[a4paper, total={6in, 9in}]{geometry}
\usepackage{hyperref}
\usepackage{setspace}
\usepackage[all,cmtip]{xy}
\usepackage{hyperref}
\usepackage{tikz}
\usepackage{tikz-cd}
\usepackage{url}

\hypersetup{
    linktoc=all,     
    linkcolor=black,  
}

\usepackage{inputenc, mathtools, amssymb, graphicx}
\usepackage{xcolor}

\newcommand{\C}{\mathbb{C}}

\newcommand{\Z}{\mathbb{Z}}

\newcommand{\R}{\mathbb{R}}

\newcommand{\eat}[1]{}
\newtheorem{theorem}{Theorem}[section]
\newtheorem{lemma}[theorem]{Lemma}
\newtheorem{prop}[theorem]{Proposition}

\newtheorem{corollary}[theorem]{Corollary}
\newtheorem{conjecture}[theorem]{Conjecture}
\newtheorem{remark}[theorem]{Remark}
\newtheorem{ex}[theorem]{Example}
\newtheorem{defn}[theorem]{Definition}

\newtheorem{assume}[theorem]{Assumption}
\newcommand{\ignore}[1]{}
\newcommand{\mf}[1]{\mathfrak{#1}}

\begin{document}
 \title{Orbit closures, stabilizer limits and intermediate $G$-varieties }
 \author{Bharat Adsul, Milind Sohoni, K. V. Subrahmanyam}
\maketitle
\begin{abstract}
In this paper we study the orbit closure problem for a reductive group $G\subseteq GL(X)$ acting on a finite dimensional vector space $V$ over $\C$. We assume that the center of $GL(X)$ lies within $G$ and acts on $V$ through a fixed non-trivial character. We study points $y,z\in V$ where (i) $z$ is obtained as the leading term of the action of a 1-parameter subgroup $\lambda (t)\subseteq G$ on $y$, and (ii) $y$ and $z$ have large distinctive stabilizers $K,H \subseteq G$. Let $O(z)$ (resp. $O(y)$) denote the $G$-orbits of $z$ (resp. $y$), and $\overline{O(z)}$ (resp. $\overline{O(y)}$) their closures, then (i) implies that  $z\in \overline{O(y)}$. 
We address the question: under what conditions can (i) and (ii) be simultaneously satisfied, i.e, there exists a 1-PS $\lambda \subseteq G$ for which $z$ is observed as a limit of $y$. 

Using $\lambda$, we develop a leading term analysis which applies to $V$ as well as to ${\cal G}= Lie(G)$ the Lie algebra of $G$ and its subalgebras ${\cal K}$ and ${\cal H}$, the Lie algebras of $K$ and $H$ respectively. 
Through this we construct the Lie algebra $\hat{\cal K} \subseteq {\cal H}$ which connects $y$ and $z$ through their Lie algebras. We develop the properties of $\hat{\cal K}$ and relate it to the action of ${\cal H}$ on $\overline{N}=V/T_z O(z)$, the normal slice to the orbit $O(z)$. 

We examine the case of {\em alignment} when a semisimple element belongs to both ${\cal H}$ and ${\cal K}$, and the conditions for the same. We illustrate some consequences of alignment and relate it to existing work in the case of the determinant and permanent. 
Next, we examine the possibility of {\em intermediate $G$-varieties} $W$ which lie between the orbit closures of $z$ and $y$, i.e. $\overline{O(z)} \subsetneq W \subsetneq O(y)$. These have a direct bearing on representation theoretic as well as geometric properties which connect $z$ and $y$.

The paper hopes to contribute to the Geometric Complexity Theory approach of addressing problems in computational complexity in theoretical computer science.

\end{abstract}

\setcounter{page}{1}
\section{Introduction}
Let $X$ be a vector space over $\C$ of dimension $n$ and let $G \subseteq GL(X)$ be a reductive algebraic group over $\C$. Furthermore, if $Z=\{ tI |t \in \C^*\}$, the center of $GL(X)$ of non-zero multiples of the identity matrix $I$, then we assume that $Z$ is a subgroup of $G$. 
Let $V$ be a finite dimensional $G$-module via a rational map $\rho : G \rightarrow GL(V)$. For a $g\in G$ and $v\in V$, let $\rho (g)\cdot v$, or simply $gv$  denote the action of $g$ on $v$ via $\rho$. We assume that $Z$ acts through a fixed non-trivial character on $V$. In other words, there is an integer $c\neq 0$ such that for any $v\in V$, we have $tI\cdot v =\rho(tI)v=t^c v$. 

Let $y\in V$ and $O(y)$ denote the $G$-orbit of $y$. Since $Z$ acts non-trivially, the closure $\overline{O(y)}$ is also a cone and its ideal $I_y \subseteq \C [V]$ is a homogeneous ideal. 
Let $\lambda : \C^* \rightarrow G$ be a one-parameter subgroup (1-PS) such that 
\begin{equation} \label{eqn1ps} \lambda (t) y = \sum_{i=d}^D t^i y_i = t^d z + t^e y_e +\ldots t^D y_D 
\end{equation}
where $z=y_d$ and $y_e$ are non-zero vectors and $d < e < \ldots < D$.  We call $z$ the leading term of $y$ under $\lambda $ and $y_e$ as the {\em tangent of approach}. 

The motivation of this paper is to study the following question. Given special elements $z, y \in V$ with large and distinctive stabilizers $H,K\subseteq G$, to determine, using just the stabilizer data, if $z$ can arise as a leading term of $y$ under a 1-PS $\lambda (t)\subseteq G$. 
By dividing by $t^d$ (which is also achieved by applying a suitable multiple of the identity), it is clear that the leading term $z$ lies within $\overline{O(y)}$, the closure of the $G$-orbit of $y\in V$. 

\subsection{Our Contributions}

Let us set up some of the background notation. 

Let ${\cal G}$ denote the Lie algebra of $G$. The central objects are ${\cal G}_y ={\cal K}$, the Lie algebra of the stabilizer $K$ of $y$ and ${\cal G}_z={\cal H}$, that of the stabilizer $H$ of $z$. Let $O(z)$ and $O(y)$ denote the $G$-orbits of $z$ and $y$ respectively. Let $\overline{O(z)}, \overline{O(y)}$ denote their closures, $I_z , I_y \subseteq \C [V]$ their ideals and $A_z , A_y$ their coordinate rings.  Let $T_z O(z) ={\cal G}\cdot z$, denote the tangent space of the orbit $O(z)$ at the point $z$. Let $\overline{N}$ be the ${\cal H}$-module $V/(T_z O(z))$ which represents a ``normal'' slice at $z$ to the tangent space $T_z O(z)$. 

Finally, the 1-PS $\lambda (t) $ allows a grading of $V$ by weights under the usual action, and ${\cal G}$ under the adjoint action. We thus have $V=\oplus_i V_i$ with $\lambda (t) v_i =t^i v_i$ for any $v_i \in V_i$. We also have ${\cal G}=\oplus_j {\cal G}_j$. For any non-zero $v=\sum_i v_i$, let the leading term $\hat{v}$ be the non-zero term $v_a$ of smallest degree. 
Similarly, for a non-zero $\mf{g}\in {\cal G}$ with $\mf{g}=\sum_j \mf{g}_j$, let $\hat{\mf{g}}$ denote the non-zero term $\mf{g}_b$ of smallest degree. Thus, in the above notation, we have $\hat{y}=z$.

We prove:
\begin{theorem}
Let $\hat{\cal K}$ be the vector space generated by $\{ \hat{\mf{k}} |\mf{k}\in {\cal K} \}$, the collection of leading terms of ${\cal K}$. Then $\hat{\cal K}$ is a Lie subalegbra of ${\cal H}$ and $dim(\hat{\cal K})=dim({\cal K})$. Moroever, $\hat{\cal K}\subseteq {\cal H}_{\overline{y_e}} $, the Lie algebra stabilizer within ${\cal H}$ of $\overline{y_e} \in \overline{N}$.  
\end{theorem}
$\hat{\cal K}$ is also the limit of $\lambda (t) {\cal K} \lambda (t)^{-1}$, the stabilizer of $y(t)=\lambda(t)y$. The above theorem brings out the role of $y_{e}$, the {\em tangent of approach} as an element of $\overline{N}$, a normal section to the orbit $O(z)$ at $z$.   
The injection $\hat{\cal K} \subseteq {\cal H}$ sets up a direct Lie algebraic connection between $y$ and its limit $z$ through their stabilizers. 

Even though ${\cal K}$ may be semisimple, $\hat{\cal K}$ demonstrates a variety of possibilities, depending on the alignment between $\lambda (t)$ and ${\cal K}$. 
Let $P(\lambda)$ be the parabolic subgroup of $G$ corresponding to $\lambda$ (see Definition~\ref{defn:Plambda}). Let $P(\lambda) = L(\lambda) U(\lambda)$ be the Levi decomposition of $P(\lambda)$. Let ${\cal P}(\lambda), {\cal L}(\lambda),{\cal U}(\lambda)$ be their Lie algebras. Finally, let $\overline{\ell} \in {\cal G}$ be the toric element such that $t^{\overline{\ell}}=\lambda (t)$. We prove:
\begin{theorem}

Let $y,z, \lambda $ and $\overline{\ell}$ be as above. Then at least one of the following holds: 
\begin{enumerate}
    \item[{\bf (A)}] Let ${\cal K}'=\hat{\cal K} \oplus \C \overline{\ell}$, then ${\cal K}' \subseteq {\cal H}$ is a Lie algebra of rank 1, i.e., the dimension of any maximal torus in ${\cal K}'$ is $1$. 
    
    or
    
    \item[{\bf (B)}] There is a semisimple element $\mf{k} \in {\cal K}$  and a (unipotent) element $u\in U(\lambda)$ such that the conjugate $\mf{k}^u \in {\cal H}$. 
\end{enumerate}
We call such an element $\mf{k}^u \in {\cal H}\cap {\cal K}^u$ as a alignment between $z$ and $y^u$. 
\end{theorem}

Alignment, i.e., the presence of a common semisimple element in ${\cal H}$ and ${\cal K}$ (or its conjugate) has important consequences for the determinant vs. permanent problem as
well as co-dimension one orbits on the boundary of the orbit of the determinant.

Let $X$ be an $n\times n$-matrix of indeterminates and $V=Sym^n (X^* )$. Consider $y=det_n (X)$ with stabilizer $K_n \subseteq GL_{n^2}$. It is well known \cite{matsushima1960espaces} that the boundary of the determinant orbit is a finite union of $G$-varieties of codimension 1, that is,  $\overline{O(y)} - O(y) = \cup_i W_i$ 
where $W_i$ is a $G$-variety of  of dimension one less than that of $\overline{O(y)}$.

\begin{corollary}
If $W_i $ equals $\overline{O(Q_i)}$ for some form $Q_i \in V$ which is obtained as a limit of a 1-PS $\lambda_i $ acting on $det_n (X)$, then either the stabilizer ${\cal H}_i$ of $Q_i$ is of rank $1$, or there is an alignment between $Q_i$ and $det_n$ (or its conjugate).
\end{corollary}

The special case of $n=3$ is analysed and the extent of alignment between $Q_i$'s and $det_3$ is illustrated. 

We also consider the case when $Y\cong \C^{m^2+1}$, is the space of $m\times m$ matrices along with an auxiliary variable $Y_{nn}$. Suppose $\phi : Y\rightarrow X$ is such that the pullback of $det_n (X)$ is the padded permanent $Y_{nn}^{n-m} perm_m (Y)$. We show that any alignment provides explicit combinatorial information on the structure of $\phi$. 

\begin{prop}
Suppose $\phi:Y\rightarrow X$ as above has an alignment then there is a {\em rectangular decomposition} ${\cal R}=\{ R_1 ,\ldots , R_r \}$ of the index set of $Y$ and ${\cal S}=\{ S_1 , \ldots , S_s \}$ of that of $X$ and a relation $\Phi \subseteq {\cal R}\times {\cal S}$ such that: 
\ $\phi (Y_{R_i})\subseteq \oplus_{S_j \in \Phi (R_i)} X_{S_j}$. 

\end{prop}

For both cases above, we show that the absence of an alignment poses exceptional conditions on $\lambda $.

The second part deals with developing the connection between ${\cal K}$ and ${\cal H}$ through {\em intermediate $G$-varieties}, defined below:
\begin{defn}
We say that the closed variety $W$ is an intermediate variety between $\overline{O(y)}$ and $\overline{O(z)}$ if $W$ is $G$-stable and $\overline{O(y)} \supseteq W \supseteq \overline{O(z)}$. We say $W$ is strict if $\overline{O(y)} \supsetneq W \supsetneq \overline{O(z)}$. 
\end{defn}

We provide two recipes to construct intermediate varieties
$\overline{O(y)} \supseteq W \supseteq \overline{O(z)}$. 

In the first case, we attempt to construct the smallest $G$-variety $W$ such that there is an $x\in W $ and a 1-PS $\mu(t)\subseteq G$ such that $z$ is the leading term of $x$ under the action of $\mu$ and $y_e$ is the tangent of approach. In other words, there is a 1-PS path lying entirely within $W$ which approaches $z$ with $y_e$ as the tangent.  

Towards this, we construct the associated graded ring for the ideal $I_z \subseteq \C [V]$ as $R=\oplus R_i$, where $R_i =I_z^i /I_z^{i+1}$. Note that $R\cong \C[V]$ as $G$-modules and $dim(R)=dim(\C[V])$ as algebras over $\C$. For an ideal $I\subseteq \C [V]$, let $\overline{I} \subseteq R$ be the graded ideal corresponding to $I$.

For any $w\in O(z)$ and $v\in T_{w} O(z)$, we define $D^k_{w,v}: I^k_z/I^{k+1}_z \rightarrow \C$ where for any $\overline{f} \in I^k_z /I^{k+1}_z$ and representative $f\in I^k_z$, $D^k_{w,v} (\overline{f})$ is the coefficient of $t^k$ in $f(w+tv)$. Thus $D^k_{w,v}$ are generalized derivations at the point $w$ in the direction $v$.  

We then have the following theorem:

\begin{theorem}
Let 
\[\overline{J}_k =\{ \overline{f} \in I^k_z /I^{k+1}_z | D^k_{gz, gy_e}(\overline{f})=0 \mbox{ for all } g\in G \} \]
Then:
\begin{enumerate}
    \item $\overline{J}_{z,y_e}=\oplus_{k\geq 1} \overline{J}_k $ is a $G$-stable ideal of $R$. Moreover, $\overline{I}_z \supseteq \overline{J}_{z,y_e} \supseteq \overline{I}_y$. 
    \item The dimension of $\overline{J}_{z,y_e}$ is $dim(G)-dim({\cal H}_{\overline{y_e}})$. 
\end{enumerate}
\end{theorem}

We know in general that $\hat{\cal K}\subseteq {\cal H}_{\overline{y_e}}$. By the above theorem, if $dim({\cal K})< dim({\cal H}_{\overline{y_e}})$, then, in the normal cone $Spec(R)$, there is indeed an intermediate variety between $\overline{O(z)}$ and $\overline{O(y)}$. 

\begin{conjecture}
If $dim({\cal K})< dim({\cal H}_{\overline{y_e}})$, then there is a strictly intermediate variety $\overline{O(z)}\subsetneq W \subsetneq \overline{O(y)}$ of dimension $dim(G)-dim({\cal H}_{\overline{y_e}})$.     
\end{conjecture}

In the second construction, we fix $\lambda $ and look at limits of elements $y' \in O(y)$ which are of the same degree as $z$. Let $\pi_d :V\rightarrow V_d$ be the projection onto the weight space. Then, we define:  
\[ Y_d =\{ y' \in O(y) \mbox{ such that $z'=\widehat{y'}$ is of degree $d$} \} \]
and $Z_d =\pi_d (Y_d)$, the set of leading terms of elements in $Y_d$. Thus $Z_d$ is the space of {\em co-limits} of $z$ and obtained from elements of $O(y)$ using $\lambda $. It is easy to see:
\begin{lemma}
Let $O(Z_d)=\{ gz' | z' \in Z_d \mbox{ and } g\in G\} $ and $\overline{O(Z_d )}$ be its closure. Then  $\overline{O(Z_d )}$ is an intermediate variety. 
\end{lemma}

The question is if the strict condition $\overline{O(z)} \subsetneq \overline{O(Z_d)}$ holds. We check this by comparing the tangent spaces of degree $d$, viz., $T_z Z_d$ and $(T_z O(z))_d={\cal G}_0 z$

\begin{defn}
Let $y, \lambda$ and $d$ be fixed as above. A element $\mf{g}\in {\cal G}$ is called a $d$-stabilizer iff $\mf{g}=\sum_{i} \mf{g}_i$ is such that $(\mf{g}y)_a =0$ for all $a<d$. Let ${\cal G}_{y,d}\subseteq {\cal G}$ be the collection of $d$-stabilizers of $y$. 
\end{defn}

In other words, ${\cal G}_{y,d}$ is the collection of ``Lie elements'' $g\in G$ such that $gy\subseteq Y_d$. 

\begin{prop} 
 Let $k=dim(H)-dim(K)$. Suppose that $y$ and $z$ are smooth within $Y_d$ and $Z_d$ respectively. Then there is a subspace $F\subseteq {\cal G}_{y,d}$ of dimension at most $k$ such that 
${\cal G}_{y,d}={\cal P}(\lambda)+{\cal K}+F$. 
Moreover, let $TW_z = \pi_d (\{  \mf{g} \cdot y | \mf{g} \in F\})$ be the leading terms of $F\cdot y$, then $T_z Z_d$ =$TW_z + {\cal G}_0 z $. 
\end{prop}

If $y$ is in the null cone of $V$ for the $G$-action and $\lambda$ is the ``optimal'' 1-PS then ${\cal G}_{y,d}={\cal P}(\lambda)$, see \cite{hesselink1979desingularizations}, Lemma 4.6. Thus $F$ measures the deviation of $\lambda $ from the optimal 1-PS which drives $y$ to $0$. If $\lambda$
 were optimal then $\overline{O(z)}=\overline{O(Z_d )}$.

\subsection{Background}

The problem of analysing the stabilizers $K$ of $y$ and $H$ of the limit $z$, arises in showing lower bounds in algebraic complexity theory, for example the permanent vs. determinant question, see \cite{mulmuley2001geometric}, Conjecture 4.3, for details. That the question of orbit closures for such special $y$ and $z$ can be settled by using purely the stabilizer data is the central thesis of Geometric Complexity Theory (GCT), see \cite{mulmuley2008geometric} and others (\cite{burgisser2011overview}, \cite{landsberg2015geometric}). 
This paper hopes to contribute to the theory by providing Lie algebraic techniques for the same. 

An earlier approach, proposed in \cite{mulmuley2001geometric} was to use the Peter-Weyl condition as follows. 
For the $H,K$ as above, the $G$-modules which appear in $A_y$ (or $A_z$) are determined by the Peter-Weyl condition, i.e., these are $G$-modules $V_{\mu}$ such that $V_{\mu}^*$ has a $K$ (resp. $H$ fixed vector). The existence of a 1-PS as above would mean that $z\in \overline{O(y)}$ and we would have a $G$-equivariant surjection $A_y \rightarrow A_z$. 
Thus, the proof of the nonexistence of a suitable 1-PS is obtained by the presence of certain $G$-modules in $A_z$ as {\em obstructions}, thereby offering a combinatorial recipe for the problem. However, for the $H,K$ in question, the mere absence or presence of certain $G$-modules as  obstructions was shown to be inadequate to prove the required exponential lower bounds, see \cite{burgisser2019no}. 
As pointed by them and other authors,  a more refined analysis of occurrences and multiplicities of $G$-modules in $A_y,A_z$, \cite{dorfler2020geometric}, \cite{ikenmeyer2020implementing} may yield the required obstructions. 

Besides GCT, there have been other algebraic approaches to lower bounds in algebraic complexity theory, and in particular to the determinant and permanent problem. 
See, for example, \cite{mignon2004quadratic}, for a quadratic lower bound which uses the curvature data for the zero sets of determinant and permanent as hypesurfaces. The more traditional approach to lower bounds problems has been to study properties of polynomials computed by circuit families with restrictions placed on their size and/or depth of these circuits. Here combinatorial and algebraic techniques are used and these approaches have met with fair success, see the recent survey by \cite{ramprasad21}. In \cite{DDS21} the authors study the closure of forms computed by algebraic ${\Sigma}^{k}\Pi\Sigma$-circuits with the top $\Sigma$-gate having constant fanin $k$. They show that forms in the closure of this circuit class have polynomial determinantal complexity, see \cite{dutta23}.  

Coming to this paper, in \cite{liemethodsGCT} we developed Lie algebraic techniques to study when $z \in \overline{O(y)}$ is obtained as the limit of a general 1-parameter family (1-PF) $\gamma (t)\subseteq G$.  
The central objects of study in \cite{liemethodsGCT} are the Lie algebras ${\cal H}$ of $H$, ${\cal K}$ of $K$ and ${\cal G}$ of $G$. We gave explicit formulas for the action of ${\cal G}$ on an appropriately chosen slice at $z$, and called this a local model at $z$. 
This was then used to construct the limiting Lie algebra $\hat{\cal K}\subseteq {\cal H}$ of ${\cal K}$ and to analyse its properties. 

The first part of this paper arose as an attempted simplification of the construction of $\hat{\cal K}$ and its properties as discussed in \cite{liemethodsGCT} when the 1-PF is actually a 1-PS. No familiarity with that paper is assumed. 

\section{Stabilizer Limits}\label{sec:leadingterm}

We recall the setting from Section 1. $G \subseteq GL(X)$ is a reductive algebraic group over $\C$ where $X$ is a vector space over $\C$ of dimension $n$. The center $Z$ of $GL(X)$ is a subgroup of $G$. $V$ is a finite dimensional $G$-module
with the center $Z$ acting as a nontrivial character on $V$. We have  $y\in V$ with stabilizer $K$ and $z\in \overline{O(y)}$ with stabilizer $H$. The ideals of $\overline{O(y)}$ and $\overline{O(z)}$ within $\C [V]$ are $I_y$ and $ I_z$ respectively. 

Let ${\cal G}=Lie(G)$ be the Lie algebra of $G$. Note that ${\cal G}$ is a subalgebra of $gl(X)=gl_n (\C)$, the Lie algebra of $n\times n$ matrices over $\C$. We recall the following basic lemma:

\begin{lemma}
For the above data, we have: 
\begin{enumerate} 
\item $G$ acts on ${\cal G}$ by conjugation and this is called the adjoint action:
$adj(g)\cdot \mf{g}= g\mf{g}g^{-1}$.
\item For the action $\rho :G\rightarrow GL(V)$, there is a Lie algebra action $\rho_1 :{\cal G} \rightarrow End(V)$ such that
for any $g\in G$ and $\mf{g}\in {\cal G}$:
\[ \rho (g) \rho_1 (\mf{g}) \rho (g)^{-1} =\rho_1 (adj(g)\cdot \mf{g})=\rho_1 (g\mf{g} g^{-1}) \]
\end{enumerate}
\end{lemma}

For any $v\in V$ and $\mf{g}\in{\cal G}$ (resp. $g\in G$), $\mf{g}v$ (resp. $gv$) will denote the action $\mf{g}$ (resp. $g$) on the element $v$ via $\rho_1 $ (resp. $\rho$). For any $v\in V$, the {\em stabilizer} of $v$ will refer to either $G_v := \{g \in G | gv = v\}$, the subgroup of $G$ fixing $v$, or to ${\cal G}_v := \{\mf{g} \in {\cal G} | \mf{g}v = 0\}$, the Lie subalgebra of ${\cal G}$ sending $v$ to zero.

\subsection{Leading term modules and algebras}

The situation we are interested in is when $y,z$ are special elements of $V$ as above and $\lambda :\C^* \rightarrow G$ is a 1-parameter subgroup (1-PS) such that:
\[ \lambda (t) y=y(t)=y_d t^d + y_e  t^e +\ldots y_D t^D \]
where $y_d =z\neq 0$. We remind the reader that in the expression on the right side above, the exponents of $t$ are increasing from left to right. The vector $z$ is the limit of $y$ under $\lambda $. 

Now $\lambda$ gives us a grading of $V$ as well as of the Lie algebra ${\cal G}$ of $G$. We state this as a lemma without proof. 

\begin{lemma}
Under the action of $\lambda$, we have a $\Z$ grading $V=\oplus_i V_i$ such that for any $v$, we have $v=\sum_i v_i$ with 
\[ \lambda (t)v =\sum_i t^i v_i \]
Similarly, we have a $\Z$-grading ${\cal G}=\oplus_j {\cal G}_j $\footnote{Note the unfortunate use of the notation ${\cal G}_y$ and ${\cal G}_j$ for the stabilizer of $y$ and also the degree $j$-component of ${\cal G}$. When $j$ is an integer it will always mean the latter.} such that for any $\mf{g}\in {\cal G}$, we have $\mf{g}=\sum_j \mf{g}_j$ and:
\[ \lambda (t) \mf{g} =\sum_j t^j \mf{g}_j  \]
Finally, if $\mf{g}v=w=\sum_i w_i$, then $w_i =\sum_{j} \mf{g}_i v_{i-j}$. 
\end{lemma}

The {\em degree} of a non-zero element $v_i \in V_i$ (resp. $\mf{g}_j \in {\cal G}_j$) is the number $i$ (resp. $j$). 

\begin{defn}
For any non-zero $v \in V$ as in the above lemma, we define the leading term as $v_a\neq 0$ of smallest degree, and denote it by $\hat{v}$. The degree $deg(v)$ is defined as the integer $a$. Similarly, for any $\mf{g}\neq 0$ as above, the leading term of $\mf{g}$ is defined as $\mf{g}_b \neq 0$ of smallest degree, and is denoted by $\hat{\mf{g}}$. Its degree is defined as $b$. As a convention, we define $deg(0_V)=deg(0_{\cal G})=\infty$. 
\end{defn}

In our motivating example $z$ is the leading term of $y$. Of importance is also the coefficient of second lowest power of $t$ in the expression, $y_e$, which we call the tangent of approach.

\begin{lemma}\label{lemma:dicho}
(A) Let $v,v' \in V$ and $deg(v) <deg(v')$ then $deg(v+v')= deg(v)$. If $deg(v)=deg(v')$ then either (i) $deg(v+v')=deg(v)$ and $\widehat{(v+v')}=\hat{v}+\hat{v'}$ or (ii) $deg(v+v')>deg(v)$. 
(B) For $\mf{g},\mf{g}'\in {\cal G}$, either (i) $deg([\mf{g}, \mf{g}'])=deg(\mf{g})+deg(\mf{g}')$ and then $[\hat{\mf{g}}, \hat{\mf{g}}']=\widehat{[\mf{g},\mf{g}']}$, or $deg([\mf{g}, \mf{g}'])>deg(\mf{g})+deg(\mf{g}')$ and then $[\hat{\mf{g}}, \hat{\mf{g}}']=0$ . 
Finally, if $\mf{g} \in {\cal G}$ and $v\in V$ are arbitrary elements then either $\hat{\mf{g}} \hat{v}=0$ or $deg(\mf{g}(v))=deg(v)+deg(\mf{g})$
\end{lemma}

\noindent 
{\bf Proof}: (A) follows from the linearity of the action of $\lambda$, i.e., $(v+v')(t)=v(t)+v'(t)$. For (B), for the first claim, note that for the adjoint action of $G$ on ${\cal G}$, we have $\lambda (t) [\mf{g}, \mf{g}'] =[ \lambda (t) \mf{g}, \lambda (t) \mf{g}']$, or in other words, $[\mf{g},\mf{g'}](t)=[\mf{g}(t), \mf{g}' (t)]$. Thus if $\mf{g}(t)=\sum_{i\geq a} \mf{g}_i t^i$ and $\mf{g}'(t)=\sum_{j\geq b} \mf{g}'_j t^j$, then we have:
\[ \begin{array}{rcl}
[\mf{g},\mf{g'}](t)&=& [\mf{g}(t), \mf{g}' (t)] \\
&=& (\mf{g}_a t^a + \ldots )(\mf{g}'_b t^b +\ldots )\\
&=& [\mf{g}_a , \mf{g}'_b ]t^{a+b}+\ldots \\ 
&=& [\hat{\mf{g}}, \hat{\mf{g}'}]t^{a+b}+\ldots 
\end{array}\]
The two cases are determined by whether $[\hat{\mf{g}}, \hat{\mf{g}'}]=0$ or not. 
For the second claim, let $deg(v)=a$ and $deg(\mf{g})=b$. We see that:
\[ \begin{array}{rcl}
\rho(\lambda(t))(\rho_1 (\mf{g})(v))&=& \rho(\lambda(t))\rho_1 (\mf{g}) \rho(\lambda (t)^{-1}) \rho (\lambda  (t))(v)) \\
&=& \mf{g}(t)v(t) \\
&=& t^{a+b} \mf{g}_b v_a + \mbox{  higher degree terms} \\
&=& t^{a+b} \hat{\mf{g}} \hat{v} + \mbox{  higher degree terms} \\
\end{array}\]
Again, which of the two conditions holds depends on whether $\hat{\mf{g}}\hat{v}=0$ or not. This proves the second claim. $\Box $

\begin{lemma}
\label{lemma:samedim}
Let $M\subseteq V$ be a subspace. Define $\hat{M}$ as the subspace generated by the set $\{ \hat{m} |m \in M \}$. Then $\hat{M}$ is a finite dimensional subspace of $V$ and $dim_{\C} (M)=dim_{\C} (\hat{M)}$. Similarly, for any subspace ${\cal L}\subseteq {\cal G}$, $\hat{\cal L}$ is a finite dimensional subspace of ${\cal G}$ and has the same dimension as ${\cal L}$.  
\end{lemma}

\noindent
{\bf Proof}: Define $M_i=\{ m\in M | deg(m)\geq i\}$, i.e., $M_i := M\cap (\oplus_{j\geq i} V_j )$. Then $M_i$ is a finite dimensional subspace of $M$ and $M_i \supseteq M_{i+1}$. Let $D=\{ deg(m)| m\in M\}$ be the set of all degrees which are seen. Let $D=\{ i_1 , \ldots ,i_k, \ldots \}$. Then $D$ is also the set of indices where $M_{i_j } \supsetneq M_{i_j +1}$. Since $M$ is finite-dimensional, we see that $D$ is finite, say $D=\{ i_1 ,\ldots ,i_k \}$. 

Let $d_j =dim(M_{i_j})-dim(M_{i_j +1})$ and $B_j =\{ m_{j,1}, \ldots, m_{j,d_j} \} \subseteq M$ be linearly independent elements of $M$ such that $M_{i_j}=M_{i_j +1}+ \C \cdot B_j$ (where $\C \cdot B_j$ is the linear space generated by the elements of $B_j $). 

Let $B=\cup_{j=1}^k B_j$. Our first claim is that (i) the elements of $B$ are linearly independent. As a contradiction, suppose that $\sum_i \alpha_i b_i =0$ for some non-zero $\alpha_i $, with $b_i \in B$. 
Let $s$ be the minimum of the degrees of $b_i$'s in this linear combination. Then restricting this linear combination to those $b_i$'s of degree $s$ gives us a non-trivial linear dependence on $M_s/M_{s+1}$. This is a contradiction to the choice of independent elements in $B$ spanning $M_s$. This proves (i).  

Next, we claim that (ii) $B$ is a basis for $M$. Suppose not, let $m\in M -\C\cdot B$ be of maximum degree among all such $m$. If $d=deg(m)$, then $d=i_j$ for some $j$. Then subtracting a suitable linear combination of the elements of $B_j$, i.e., $m'=m-\sum_r \alpha_r m_{i_j ,r}$ will give us an element $m'$ of a higher degree and an element outside $\C \cdot B$. But this contradicts the choice of $m$. This proves (ii). 

Finally, (iii) we claim that $\hat{B}$ is a basis for $\hat{M}$. Clearly, all leading terms of $\hat{M}$ have degrees from the set $D$. For any $\hat{m}$ of degree $d=i_j$, we know that there is an element $m''=\sum_r \alpha_r m_{i_j ,r}$ of $M$ such that $m'=m-m''$ is either of a higher degree or is zero. 
In either case, we have $\hat{m}=\sum_r \alpha_r \widehat{m_{i_j ,r}}$. 
Coming to the linear independence of $\hat{B}$, note that $\hat{B}=\dot{\cup}_j \hat{B}_j$, where each $\hat{B}_j $ is a subset of different graded components $V_j$. Thus any linear dependence must be purely within $\hat{B}_j$. But that would force a linear dependence on $B_j$. This proves (iii) and the assertion about $\hat{M}$. The $\hat{\cal L}$ assertion is similarly proved. $\Box $

With this definition, we have the following lemmas.

\begin{lemma}\label{lemma:general-1}
Let ${\cal K}$ be a Lie subalgebra of ${\cal G}$ and $N\subseteq V$ a ${\cal K}$-module. Then (i) $\hat{\cal K}$ is a graded Lie subalgebra of ${\cal G}$, and $dim_{\C} (\hat{\cal K})=dim_{\C} ({\cal K})$, (ii) $\hat{N}\subseteq V$ is a $\hat{\cal K}$-module with $dim_{\C} \hat{N}=dim_{\C} N$. 
\end{lemma}

{\bf Proof}: Let us first prove that $\hat{\cal K}$ is a Lie subalgebra. For that, it is adequate to show that $[\mf{k}_1 ,\mf{k}_2] \in \hat{\cal K}$ for any leading terms $\mf{k}_1 ,\mf{k}_2 \in \hat{\cal K}$. Given such elements, there are elements $\mf{k}'_1, \mf{k}'_2 \in {\cal K}$ such that $\mf{k}_i = \hat{\mf{k}'_i}$. 
Let the degrees of the leading terms be $d_1$ and $d_2$. Consider the element $\mf{k}'=[\mf{k}'_1 ,\mf{k'}_2] \in {\cal K}$. Note that $\mf{k}'(t)=[\mf{k'_1} (t),\mf{k}'_2 (t)]$, and hence, by Lemma \ref{lemma:dicho}, 
either (a) the leading term of $\mf{k}'$ is of degree $d_1 +d_2$, in which case, we have $\mf{k}=[\mf{k}_1 ,\mf{k}_2]$ is the leading term of $\mf{k}'$, and so is an element of $\hat{\cal K}$, or (b), the leading term is of a higher degree, in which case $[\mf{k}_1 ,\mf{k}_2]=0$. 
This proves that $\hat{\cal K}$ is a Lie subalgebra. That $\hat{\cal K}$ is graded is clear since it is generated by leading terms, which are homogeneous. Finally, $dim_{\C} (\hat{\cal K})=dim_{\C } ({\cal K})$ follows from lemma \ref{lemma:samedim}. 

Let us now prove that $\hat{N}$ is a $\hat{\cal K}$-module. For that, take any leading term $\mf{k}\in \hat{\cal K}$ and $n\in \hat{N}$. Let $\mf{k}'\in {\cal K}$ and $n' \in N$ be such that $\hat{\mf{k'}}=\mf{k}$ and $\hat{n'}=n$. 
We see that $(\mf{k}'n')(t)= \mf{k}'(t) n'(t)$, whence again, either $\mf{k}n=0$ or it equals $\widehat{(\mf{k}'n')} \in \hat{N}$. Again, by lemma \ref{lemma:samedim}, $dim_{\C} (N)=dim_{\C} (\hat{N})$ . This proves the lemma. 
$\Box $

\begin{lemma} \label{lemma:basicm-2}
Let $v$ be an arbitrary element of $V$ and  $\mf{g}\in {\cal G}$ is such that $\mf{g}\cdot v=0$. Then $\hat{\mf{g}}\cdot \hat{v}=0$. In other words $\widehat{({\cal G}_v )}\subseteq {\cal G}_{\hat{v}}$.   
\end{lemma}

\noindent
{\bf Proof}: Suppose that $v$ is of degree $a$ and $v_a =\hat{v}$. Similarly, suppose that $\mf{g}\in {\cal G}_v$, the stabilizer of $v$, is of degree $b$ and $\mf{g}_b =\hat{\mf{g}}$. Then $\mf{g}v=0$ implies $\mf{g}(t)v(t)=0$ as well. 
This implies that terms of all degree in the product are zero, in particular, that of degree $a+b$, viz. $\mf{g}_a v_b =0$ and thus $\hat{\mf{g}}\in {\cal G}_{\hat{v}}$. This proves the assertion. $\Box $

\begin{remark} Let $Mod_{\cal K}(V)$ (resp. $Mod_{\hat{\cal K}} (V)$) be the collection of ${\cal K}$ (resp. $\hat{\cal K}$) submodules of $W$. Then lemmas \ref{lemma:general-1} and lemma \ref{lemma:basicm-2} set up a map $Mod_{\cal K} (V) \stackrel{\lambda}{\rightarrow} Mod_{\hat{\cal K}}(V)$, where the ${\cal K}$-module $N$ goes to the $\hat{\cal K}$-module  $\hat{N}$, which is of the same dimension. Moreover, if $N$ has a ${\cal K}$-fixed point then $\hat{N}$ has a $\hat{\cal K}$-fixed point. 
\end{remark}

We illustrate Lemma \ref{lemma:general-1} with two examples. The first example illustrates the
computation of the leading term algebra $\hat{\cal K}$ from ${\cal K}$ and its dependence of $\hat{{\cal K}}$ on the {\it alignment} of ${\cal K}$ with respect to the 1-parameter subgroup.

\begin{ex} \label{ex:hatkchange}
Let us consider $G=GL_4 (\C)$ and the 1-parameter subgroup $\lambda (t)$ below. The action of $\lambda $ on a typical element $M \in {\cal G}$ is given below, where each $M_{ij}$ is a $2 \times 2$-matrix. Note that the degrees which occur are $-1,0$
and $1$, thus we have the spaces ${\cal G}_{-1}, {\cal G}_{0}$ and ${\cal G}_1$, of matrices with leading terms of degree $-1, 0$ and $1$ respectively. 

\[ 
\lambda(t)=
\left[ \begin{array}{cccc}
1 & 0 & 0 & 0\\
0 & 1 & 0 & 0\\
0 & 0 & t & 0 \\
0 & 0 & 0 & t \end{array} \right], 
\lambda (t) M \lambda(t)^{-1} = \left[ \begin{array}{cc}
M_{11} & t^{-1} M_{12} \\
tM_{21} & M_{22} 
\end{array} \right]
\]
Let ${\cal K}$ be as shown below and let us construct $\hat{\cal K}$. Let ${\cal K}_i ={\cal G}_{j\geq i} \cap {\cal K}$. We then have $dim({\cal K}_{-1})=4, dim({\cal K}_0 )=4$ and $dim({\cal K}_1 )=0$. Thus, we get a basis $B=B_0$ of ${\cal K}$ of dimension $4$. Moreover, since each element of $B_0$ is homogeneous, we have $\hat{B_0}=B_0$ and $\hat{\cal K}={\cal K}$. 

\[ {\cal K}= \left[ \begin{array}{cccc}
a & b & 0 & 0 \\
c & d & 0 & 0\\
0 & 0  & a & b \\
0 & 0 & c & d 
\end{array} \right]
\: \: 
\lambda(t)=
\left[ \begin{array}{cccc}
1 & 0 & 0 & 0\\
0 & 1 & 0 & 0\\
0 & 0 & t & 0 \\
0 & 0 & 0 & t \end{array} \right] \]

The leading term Lie algebra $\hat{{\cal K}}$ of ${\cal K}$ depends intimately on $\lambda (t)$ and may change dramatically under conjugation. Let $A$ be as shown below and ${\cal K}'=A{\cal K}A^{-1}$. This gives us ${\cal K}'$ as shown below. Let us compute $\hat{\cal K}'$. 
\[ A=\left[ \begin{array}{cccc}
1 & 0 & 0 & 1\\
0 & 1 & 0 & 0\\
0 & 0 & 1 & 0 \\
0 & 0 & 0 & 1 \end{array} \right] \: \: 
{\cal K}'= \left[ \begin{array}{cccc}
a & b & c & d-a \\
c & d & 0 & -c\\
0 & 0  & a & b \\
0 & 0 & c & d 
\end{array} \right] \;\;
\]

As before, let ${\cal K}'_i = {\cal K}'\cap {\cal G}_{j\geq i}$ and note that $dim({\cal K}'_{-1})=4$. Now ${\cal K}'_0$ consists of all elements $\mf{k}\in {\cal K}'$ which have no term of degree $-1$. This forces $d=a$ and $c=0$, making $dim({\cal K}'_0 )=2$.Finally $dim({\cal K}'_1)=0$. Thus $\hat{\cal K}'$ is what is given below, with $r,s,t,u \in \C$. Note that while $\hat{\cal K}$ is reductive, $\hat{\cal K}'$ is a solvable Lie algebra. 
\[ \hat{\cal K}'= \left[ \begin{array}{cccc}
u & t & s & r \\
0 & u & 0 & -s\\
0 & 0  & u & t \\
0 & 0 & 0 & u 
\end{array} \right]
\]
\end{ex}
\begin{ex} \label{ex:O3O2}
This example illustrates the consequences of Lemma \ref{lemma:basicm-2}. Consider a set of indeterminates $\{ x,y,z\}$ and $X=\C \cdot \{x,y,z\}$. Let $GL_3$ acts on $X$ in the natural way, and let $V=Sym^4 (X)$. Consider  $f=(x^2 +y^2 +z^2)^2 \in V$. The stabilizer algebra ${\cal G}_f$ is $3$-dimensional and is given below. Consider next the 1-PS $\lambda (t) \subseteq GL(X)$ given by $\lambda (x)=x, \lambda (y)=y$ and $\lambda (z)=tz$, as shown below. We have $g=\hat{f}$ is the leading term of $((x^2+y^2+t^2 z^2)^2)$ and is $(x^2 +y^2)^2$. The stabilizer ${\cal G}_g$ is $4$-dimensional and is shown below (with $a,b,c,d \in \C$). Note that the last column of ${\cal G}_g$ is given by the operator $\nabla = (bx+cy+dz)\frac{\partial}{\partial z}$ which certainly stabilizes $g$. 
\[ {\cal G}_f= \left[ \begin{array}{ccc}
0 & a & b \\
-a & 0 & c \\
-b & -c & 0 \end{array} \right] 
\: \: \lambda(t)=
\left[ \begin{array}{ccc}
1 & 0 & 0 \\
0 & 1 & 0 \\
0 & 0 & t \end{array} \right] \: \: 
{\cal G}_g= \left[ \begin{array}{ccc}
0 & a & b \\
-a & 0 & c\\
0 & 0  & d \end{array} \right] 
\]
Let us compute $\widehat{{\cal G}_f}$ as leading terms of  ${\cal G}_f (t)$ below: 
\[ \lambda(t) {\cal G}_f \lambda (t)^{-1}= \left[ \begin{array}{rrr}
0 & a & t^{-1} b \\
-a & 0 & t^{-1} c \\
-tb & -tc & 0 \end{array} \right]
\]
Reasoning as we did in the previous example, $\widehat{{\cal G}_f}$is the $3$-dimensional Lie algebra of matrices (with $a,b,c\in \C$) shown below. 
\[ \left[ \begin{array}{rcc}
0 & a & b \\
-a & 0 & c \\
0 & 0 & 0 \end{array} \right] \subseteq {\cal G}_g
\]
Thus the "leading term" operation applied to ${\cal G}_f$, the stabilizer of $f$, inserts $\widehat{{\cal G}_f}$ into the stabilizer ${\cal G}_g$ of the limit $g$.  
\end{ex}

\begin{prop}
\label{prop:stronger}
Recall that $K=G_y$ and
$H=G_z$ are stabilizer groups of $y$ and $z$ resp.~and ${\cal K}={\cal G}_y$ and ${\cal H}={\cal G}_z$,  their stabilizer Lie algebras. 
\begin{enumerate}
    
    \item Let $T_z O(z) ={\cal G}\cdot z$, be the tangent space of the orbit $O(z)=G\cdot z$ at the point $z$. Then $V/(T_z O(z))$ is an ${\cal H}$-module. We call this the $\star$-action. 
\item Let ${\cal H}_{\overline{y_e}}$ be the stabilizer of $\overline{y_e} \in V/T_z (O(z))$ for the above action. Then the subalgebra $\hat{\cal K} \subseteq {\cal H}_{\overline{y_e}}$, the stabilizer of the image of the tangent of approach in $V/T_z (O(z))$.
\end{enumerate}
\end{prop}

\noindent 
{\bf Proof}: Since ${\cal H}={\cal G}_z$, it fixes the tangent space $T_z (O(z))$. The space $V/(T_z (O(z))$ is the quotient of ${\cal H}$-modules and hence is itself an ${\cal H}$-module. Thus (1) is clear. 

Next, let $\mf{k}\in {\cal K}$ be arbitrary. Then $\mf{k}y=0$ implies $\mf{k}(t)y(t)=0$. If $\mf{k}(t)=\mf{k}_a t^a +\mf{k}_{a+1}+\ldots$ and $y(t)= y_d t^d +y_e t^e +\ldots $, then we have:
\[ (\mf{k}_a t^a +\mf{k}_{a+1} t^{a+1}+\ldots )(y_d t^d +y_e t^e +\ldots )=0 \]
Examining the terms of degree $a+d, a+d+1,\ldots,a+e-1 , a+e$ in this product, we have:
\[ \begin{array}{rl} 
\mf{k}_{a+i} y_d=0 & \mbox{ for $i=0,\ldots e-d-1$}\\ 
\mf{k}_a y_e + \mf{k}_{a+e-d} y_d =0& \end{array} \] 
This tells us that $\mf{k}_{a+i} \in {\cal H}$ for $i=0,\ldots, e-d-1$. Since $\mf{k}_a \in {\cal H}$ and $\mf{k}_a y_e \in {\cal G} y_d$, we have $ \mf{k}_a \cdot \overline{y_{e}}=0$ and thus $\mf{k}_a =\hat{\mf{k}} \in {\cal H}_{\overline{y_e}}$. 
Since elements of the type $\hat{\mf{k}}$, with $\mf{k}\in {\cal K}$ generate $\hat{\cal K}$, we have proved (2). $\Box $
\begin{ex}
We continue with example ~\ref{ex:O3O2} to illustrate the proposition just proved.
The limit of $(x^2 + y^2 + z^2)^2$ under the given 1-PS is $g = (x^2 + y^2)^2$ and the tangent of approach is $y_e = (x^2 + y^2)z^2$. The tangent space to the orbit of $g$ contains the form $(x^2 + y^2)(2bxz + 2cyz)$ for arbitrary $b,c \in \C$. This can be seen by applying the following differential operator to $(x^2 + y^2)^2$.
\[ \left[ \begin{array}{rcc}
0 & 0 & 0 \\
0 & 0 & 0 \\
b & c & 0 \end{array} \right] 
\]
A generic  element of $\widehat{{\cal G}_f}$ corresponds to the differential form 
\[ ay \frac{\partial}{\partial x} - ax \frac{\partial}{\partial y} + bx \frac{\partial}{\partial z} + cy \frac{\partial}{\partial z}.\]
Applying this to $y_e$ gives us
$(x^2+y^2)(2bxz + 2cyz)$. Now this is zero in $V/T_g(O(g))$, since $(x^2+y^2)(2bxz + 2cyz) \in T_g(O(g))$ as we have just seen. So $\overline{y_e}$ is stabilized by 
$\widehat{{\cal G}_f}$. In fact, the full stabilizer of $\overline{y_e}$ is 
$\widehat{{\cal G}_f}$.
\end{ex}
\begin{remark}
We use the notation of the motivating example - $y \in V$ picks up $z$ as the limit of a 1-PS $\lambda$. Let $N$ be a complement to $T_z O(z)$ within $V$. We denote by $\overline{N}$, the $H$-module (and therefore ${\cal H}$-module) $V/T_z O(z)$.
\end{remark}
\subsection{Alignment} \label{subsec:1PS}

In this subsection, we look at elements which are in both ${\cal H}$ and ${\cal K}$ or in ${\cal H}$ and a conjugate of ${\cal K}$ . Such common elements will indicate common nested subspaces in the two stabilizers. We define:
\begin{defn}
    For $y,z$ and $\lambda$ as above, a semisimple element $\mf{k}\in {\cal H}\cap {\cal K}^g$, for some $g\in G$ is called an alignment. 
\end{defn}

Recall that the center $Z=\{ tI|t\in\C^*\}$ acts non-trivially on $V$. Whence, for our 1-PS $\lambda (t)$, there is a 1-PS $\lambda '(t)=t^a \lambda(t)$ such that: 
\[ \lambda' (t)y=y_d + y_{e} t^{e-d} +\ldots +y_{D} t^{D-d}\] 
with $y_d =z$ as before. Thus $z$ is stabilized by $\lambda' (t)$ and thus $\lambda'(t)\subseteq H$.  

We begin with a few definitions.

\begin{defn} \label{defn:ell}
Define $\overline{\ell} \in {\cal G}$ as the element such that $t^{\overline{\ell}} =\lambda' (t)$. 
\end{defn}

Note that $\overline{\ell} \cdot v_i=(i-d) v_i$ for all $v_i \in V_i$. 
\begin{defn}
\label{defn:Plambda}
Let $P(\lambda )$ be defined as below:
\[ P(\lambda )=\{ g \in G | \lim_{t\rightarrow 0} \lambda(t) g \lambda (t)^{-1} \mbox{ exists} \} \]
and let $U(\lambda)$ be its unipotent radical. Let $L(\lambda)$ be elements of $G$ which commute with $\lambda(t)$. Then $L(\lambda )$ is also a specially identified reductive complement to $U(\lambda)$. Let ${\cal P}(\lambda), {\cal L}(\lambda)$ and ${\cal U}(\lambda)$ be the Lie algebras of $P(\lambda ), L(\lambda)$ and $U(\lambda)$. \footnote{Note that $P(\lambda)=P(\lambda')$ and so on.} 
\end{defn}
Note that ${\cal P}(\lambda)=\oplus_{a\geq 0} {\cal G}_a$, ${\cal L}(\lambda )={\cal G}_0$ and ${\cal U}(\lambda )=\oplus_{a>0} {\cal G}_a$,\cite[Section 2]{kempf1978instability}.

We now begin with examining further the connection between $K$ and $H$ through ${\cal K}$ and $\hat{\cal K}$. More specifically, we examine which elements of $K$ (or its conjugate) descend into $H$. 

\begin{lemma}\label{lemma:pure}
With $y,z,\lambda$ as above, (i) if $R\subseteq K$ commutes with $\lambda$, then $R\subseteq H$. Furthermore,
$Lie(R) \subseteq \hat{\cal K} \subseteq {\cal H}_{\overline{y_e}} \subseteq {\cal H}$, (ii) if $\mf{g} \in {\cal L}(\lambda) \cap {\cal K}$ then $\mf{g} \in \hat{\cal K} \subseteq {\cal H}_{\overline{y_e}} \subseteq {\cal H}$. 
\end{lemma}

\noindent 
{\bf Proof}: If $\sigma \in K$ commutes with $\lambda$, then we have $\sigma \lambda (t) y= \lambda (t)\sigma y=\lambda (t) y$. Whence $\sigma $ stabilizes each degree component of $\lambda (t)y$, and therefore $z$ and $y_e$ as well. Part (i) follows. For (ii) note that, by its very definition, every element of ${\cal L}(\lambda)$ commutes with $\lambda $. $\Box $
\begin{prop} \label{prop:unistablegeneral}
Let $y,z, \lambda $ be as above. Suppose that $\mf{k} \in {\cal P}(\lambda) \cap {\cal K}$ is a semi-simple element, then there is a unipotent element $u\in U(\lambda)$ such that:
\begin{enumerate}
\item $\mf{k}^u =u\mf{k}u^{-1} \in {\cal L} (\lambda)$, i.e., it commutes with $\overline{\ell}$ and $\lambda$, and it stabilizes $y^u =u\cdot y$. 
\item Moreover $z$ is the leading term of $y^u$, i.e.,
    \begin{equation} \label{eqn:apply}        \lambda (t) y^u = t^{d} z + \mbox{ higher degree terms }
    \end{equation}
    Thus, $z=\hat{y^u}$ and $\mf{k}^u \in {\cal H} \cap {\cal G}_{y^u}$. 
    \end{enumerate}
In other words, $\mf{k}^u$ is an alignment between $z$ and $y^u$. 
\end{prop}

\noindent 
{\bf Proof}: 
Note that $P(\lambda )=L(\lambda ) U(\lambda )=U(\lambda ) L(\lambda )$ is a Levi factorization, with $L(\lambda )$ as a reductive complement. Since $\mf{k}$ is a semisimple element of ${\cal P}(\lambda)$ there is a maximal torus $T'$ of ${\cal P}(\lambda)$ containing $\mf{k}$. 
Since maximal tori in $P(\lambda)$ are $U(\lambda)$-conjugate, $T'$ is conjugate to the maximal torus in $L(\lambda)$.  So there is a $u\in U(\lambda)$ so that $\mf{k}^u=u\mf{k} u^{-1} \in {\cal L}(\lambda)$. 
Moreover, it is straightforward that $\mf{k}^u$ stabilizes $y^u =u\cdot y$. This proves (1). 

Suppose $y$ has the weight decomposition:
\begin{equation} 
\label{eqn:apply2}       
y = z + y_e + w
    \end{equation}
where $y_e \in V_e$ and $w \in \oplus_{i>e} V_i$. On applying $u$ to Eq. \ref{eqn:apply2} we have:

\[ u\cdot y = u\cdot z+u\cdot y_e +uw \]
Since $u$ is unipotent, for all $j$, we have  $u\cdot (\oplus_{i\geq j} V_i )\subseteq \oplus_{i\geq j} V_i $. 
Thus we have the graded expression:
\[ u\cdot y= z + y'_{d+1} +  w'\]
with $y'_{d+1} \in V_{d+1}$ and $w'\in \oplus_{i>d+1} V_i $. Thus, $z$ continues to be the leading term of $y^u$. In other words, we have:
\[ \lambda (t) (y^u )= t^{d} z+ t^{d+1} y'_{d+1} + \mbox{higher degree terms} \]
Now $\mf{k}^u$ commutes with $\lambda$, $\mf{k}^u$ stabilizes $y^u$ and $z=\hat{y^u}$ with stabilizer ${\cal H}$. It follows from the previous lemma that $\mf{k}^u \in \widehat{{\cal G}_{y^u}}\subseteq {\cal H}$. Thus $\mf{k}^u $ is the required alignment.
$\Box$ 

Thus, conjugates of semisimple elements $\mf{k} \in {\cal K}\cap {\cal P}(\lambda)$ are also elements of ${\cal H}$. Note that $O(y^u)=O(y)$. The next proposition handles the general case about the intersection ${\cal P}(\lambda ) \cap {\cal K}$. 

\begin{remark} \label{remark:generic}
Let ${\cal G}_- =\oplus_{i<0} {\cal G}_i$, be the complement of ${\cal P}(\lambda)$ and $\Pi_{-}:{\cal G} \rightarrow  {\cal G}_-$ be the projection. 
The condition that $\hat{\cal K} \cap {\cal P}(\lambda ) \neq 0$ is tantamount to saying that $dim(\Pi_{-} ({\cal K}))< dim({\cal K})$, i.e., $\lambda $ is not {\em generically placed} with respect to ${\cal K}$. 
\end{remark}

\begin{prop} \label{prop:stabgeneral}
Let $y,z, \lambda $ and $\overline{\ell}$ be as above. Then at least one of the following holds: 
\begin{enumerate}
    \item[{\bf (A)}] Let ${\cal K}'=\hat{\cal K} \oplus \C \overline{\ell}$, then ${\cal K}'$ is a Lie algebra of rank 1, i.e., the dimension of any maximal torus in ${\cal K}'$ is $1$. 
    
    or
    
    \item[{\bf (B)}] there is a unipotent element $u\in U(\lambda)$ and a semisimple element $\mf{k} \in {\cal K}^u$ such that $\mf{k} \in {\cal H}$. In other words, there is an alignment between $z$ and $y^u$. 
\end{enumerate}

\end{prop}

\noindent 
{\bf Proof}: Note that ${\cal K}'=\hat{\cal K} \oplus \C \overline{\ell}$ is indeed a Lie algebra since $[\overline{\ell} ,\hat{\cal K}] \subseteq \hat{\cal K}$. 

\noindent 
{\bf Case 1}: 
Suppose that $dim(\Pi_{-} ({\cal K})) =dim({\cal K})$. If this happens, then $\hat{\cal K} \subseteq {\cal G}_-$ and $\hat{\cal K}$ is nilpotent. Even more, ${\cal K}' =\hat{\cal K} \oplus \C \overline{\ell}$ is a Levi decomposition of ${\cal K}'$. Thus (A) holds. 

\noindent 
{\bf Case 2}: On the other hand, if 
$dim(\Pi_{-} ({\cal K})) <dim({\cal K})$, then ${\cal K}\cap {\cal P}(\lambda) \neq 0$.

\begin{enumerate}
\item[{\bf 2a}]
There is a semisimple element $\mf{k} \in {\cal P}(\lambda )\cap {\cal K}$.  Then, by Prop. \ref{prop:unistablegeneral}, there is indeed a $u\in U(\lambda)$ such that $\mf{k}^u \in {\cal H}$, and (B) hold. 

\item[{\bf 2b}] We are left with the case that there are no semisimple elements in ${\cal P}(\lambda) \cap {\cal K}$. Suppose now that the rank of $\hat{\cal K} \oplus \overline{\ell}>1$. Then, there is an element $\mf{k}\in {\cal K}$ such that $\hat{\mf{k}}$ is semisimple and $[\overline{\ell} ,\hat{\mf{k}}]=0$. This implies that the leading term of $\mf{k}$ is of degree zero and therefore $\mf{k}\in {\cal P}(\lambda) \cap {\cal K}$. Thus $\mf{k}=\hat{\mf{k}}+ \mf{k}_+$, where $\mf{k}_+ \in {\cal U}(\lambda)$. 

Note that ${\cal P}(\lambda) \cap {\cal K}$ is the Lie algebra of an algebraic subgroup of $GL(X)$. So the Jordon-Chevalley decomposition holds for all elements $\mf{k} \in {\cal P}(\lambda )\cap {\cal K}$, see \cite[Chapter 1, Section 4]{borel1991linear}. Hence, $\mf{k}$ may be written uniquely as a sum of a semisimple element $\mf{k}_s$ and a nilpotent element $\mf{k}_n$, $\mf{k} = \mf{k}_s + \mf{k}_n$, with $[\mf{k}_s, \mf{k}_n]=0$, and $\mf{k}_s, \mf{k}_n \in {\cal P}(\lambda)\cap {\cal K}$.  Since there are no semisimple elements in ${\cal P}(\lambda)\cap {\cal K}$, we must have $\mf{k}=\mf{k}_n$. Since $\mf{k}$ is now a nilpotent matrix, we have $\mf{k}^k=0$ for some $k>0$. This implies that $\hat{\mf{k}}^k=0$ as well. That contradicts the assumption that $\hat{\mf{k}}$ is semisimple. 

\end{enumerate}
This proves the proposition. $\Box $

\subsection{Alignment: The boundary of the general determinant $det_n$ and $det_3$}
\label{subsection:det3}
In this and the next section we study forms of interest in GCT and ask whether the results developed connect with existing work. We show that this is the case when there is alignment. We separately address the case when there is no alignment.

Let $X=(X_{ij})$ be an $n\times n$-matrix of indeterminates and $V=Sym^n (X^*)$. Consider $y=det_n (X)$ with stabilizer $K_n \subseteq GL_{n^2}$. It is well known that since the stabilizer of $det_n$ is reductive, $\overline{O(y)} - O(y) = \cup_i W_i$, is the union of closed $G$-varieties of co-dimension 1, , see for example  \cite[4.2]{burgisser2011overview}.

We then have the corollary:

\begin{corollary}
Suppose that $W_i =\overline{O(Q_i)}$, for a form $Q_i \in \overline{O(det_n )}$ obtained as a limit of a 1-PS $\lambda_i $, then either the stabilizer ${\cal H}_i$ of $Q_i$ is of rank $1$, or there is an alignment between $Q_i$ and a conjugate of $det_n$.    
\end{corollary}

\noindent 
{\bf Proof}: let ${\cal H}_i$ be the stabilizer of $Q_i$. Then since $dim({\cal H}_i) =dim({\cal K}_n )+1$, we have ${\cal H}_i ={\cal K}'= \hat{\cal K}_n \oplus \C \overline{\ell}$. By Prop. \ref{prop:stabgeneral}, either ${\cal H}_i$ must be of rank 1 or there must be a semisimple $\mf{k}\in {\cal K}_n$ whose conjugate $\mf{k}^u \in {\cal H}_i$. This proves the assertion. $\Box$

\begin{corollary}
    Suppose that there is indeed an alignment, then ${\cal R}=({\cal H}_i)_0 \cap {\cal K}_n^u \neq 0$, and thus there is a common subgroup $R\subseteq H_i\cap K_n^u$ of atleast rank $1$ which stabilizes both $Q$ as well as $(det_n)^u$.  
\end{corollary}

\begin{remark} \label{remark:recipe}
The above result provides a recipe for constructing and testing possible boundary forms $Q$s which are aligned with $det_n$. The steps are:
\begin{enumerate}
    \item Pick a subgroup $R\subseteq K_n $ which decomposes $X=\oplus_{i=1}^r X_r$. 
    \item Pick a coarsening  $X=\oplus_{i=1}^s Y_i$ of the partition above. Pick $e_1 ,\ldots , e_s$ suitably and construct $\lambda (t)$ such that $\lambda (t)(Y_i )=t^{e_i} Y_i$. 
    \item Compute the leading term $Q=\widehat{det_n}^{\lambda}$. Compute the dimension of ${\cal G}_Q$. If this equals $dim({\cal K}_n)+1$, then ${\cal G}_Q=\hat{\cal K}_n \oplus \C \ell$ and $Q$ is a boundary form.
\end{enumerate}
\end{remark}

We illustrate the above recipe for the $3\times 3$-determinant. These calculations are inspired by the work in \cite{huttenhain2016boundary}. 

Let $X=\{ x_1 ,\ldots , x_9 \}$ and $V=Sym^3 (X)$ be the space of homogeneous forms of degree $3$ acted upon by $GL(X)$. The $3\times 3$-determinant, $det_3 (X)$ is a special element as given below:
\[ det_3 (X)=det\left( \left[
\begin{array}{ccc}
x_1 & x_2 & x_3 \\
x_4 & x_5 & x_6 \\
x_7 & x_8 & x_9 \\ \end{array} \right] \right) 
\]

The stabilizer $K_3 \subseteq GL(X)$ of $det_3 (X)$ is given by transformations (i) $X\rightarrow AXB^{-1}$, where $A,B\in GL_3$ with $det(AB^{-1})=1$, and (ii) $X\rightarrow X^T$. The dimension of $K$ is $16$. Sitting within $K_3$ are two groups $R_1 =\{ X\rightarrow AXA^{-1}| A \in GL_3\}$ and $R_2=\{ X\rightarrow AXA^T | A\in SL_3 \}$. 

H\"uttenhain and Lairez \cite{huttenhain2016boundary} have proved that the boundary of the $GL(9)$-orbit of $det_3$ has two irreducible components. Moreover, these are the orbit closures of two forms $Q_1$ and $Q_2$ given below. 

\[ \begin{array}{rcl}
Q_1 (X)&=&det\left( \left[
\begin{array}{ccc}
x_1 & x_2 & x_3 \\
x_4 & x_5 & x_6 \\
x_7 & x_8 & -x_5 -x_1  \\ \end{array} \right] \right) \\ 
Q_2 (X) &=& 2(x_4 x_1^2 +x_5 x_2^2 +x_6 x_3^2 + x_7 x_1 x_2 +x_8 x_2 x_3 +x_9 x_1 x_3) \end{array}
\]

It is easy to check that $Q_1$ and $Q_2$ arise from $R_1$ and $R_2$ using Remark \ref{remark:recipe}.

\begin{prop} \label{prop:det3}
With above notation, we have:
\begin{enumerate}
    \item The space $X$ under $R_1$ decomposes as $X=X_0 \oplus cI$, the space of trace-zero matrices, and multiples of the identity. Under $R_2$, we have $X=X_a \oplus X_s $, where $X_a$ is the space of antisymmetric matrices and $X_s$, that of symmetric matrices.
    \item Let $\lambda_1$ be such that $\lambda_1 (x_0)=x_0$ for all $x_0 \in X_0$, while $\lambda_1 (I)=tI$. Similarly, let $\lambda_2 (x_a)=x_a$ and $\lambda_2 (x_s)=tx_s$ for all $x_a \in X_a$ and $x_s \in X_s$. Then $\widehat{det_3}^1=Q_1$ and $\widehat{det_3}^2=Q_2$ upto conjugates.
\end{enumerate}
    
\end{prop}

Let us now describe the stabilizers of the forms $Q_1$ and $Q_2$ above for $n=3$. 

\begin{lemma}\label{lemma:det3}
The stabilizer of ${\cal H}_i$ of $Q_i$ within $gl(X)$ has dimension 17. Moreover, it may be expressed as ${\cal H}_i = \ell^i \oplus \widehat{\cal K}^i$, where $t^{\ell^i} =\lambda_i (t)$ and $\widehat{\cal K}^i$ is the leading term algebra of ${\cal K}$ under $\lambda_i (t)$. Moreover, $({\cal H}_i){\overline{y_e}}$, the stabilizer of the tangent of approach equals $ \widehat{\cal K}^i$. Finally, $Lie(R_i) \subseteq {\cal H}_i$. 
\end{lemma}

The proofs of Propositon~\ref{prop:det3} and Lemma~\ref{lemma:det3} are computations. 

What if $Q$ has no alignment? Is $\lambda $ special in this case too?  This is partly answered by the following proposition. 

\begin{prop} \label{prop:nsquare}
Let $T\subseteq GL(X)$ be a maximal torus containing $\lambda$. Let $det_n =\sum_\alpha a_{\alpha} X^{\alpha}$ be the expression for the determinant in this basis. For any monomial index $\alpha$, let  $\xi (\alpha)$ denote its $T$-weight. Define $\Xi_T (det_n )=\{ \xi (\alpha) | a_{\alpha \neq 0} \}$ as the support of $det_n$ for this $T$. Similarly, define $\Xi_T (Q)$ as the support of the leading term $Q$. 
Then, in the absence of an alignment (i) the dimension of the $\R$-vector space formed by $\Xi_T (det_n )$ is $n^2$ while that formed by $\Xi_T (Q)$ is $n^2-1$. Moreover, $\langle \overline{\ell} , \chi \rangle \geq 0$ for all $\chi \in \Xi_T (det_n)$ but $\langle \overline{\ell}, \beta \rangle =0$ for all $\beta \in \Xi_T (Q)$.  
\end{prop}
It is easy to see that if $T\cap K_n \neq \{ I \}$, the identity, then  $dim(\Xi_{T} (det_n))<n^2$.  See Section \ref{subsection:character} for details.

If a boundary form $Q$ is obtained by a suitable $\lambda $ then a conjugate $\lambda' $ is available in any maximal torus $T\subseteq GL(X)$. Then, the requirement of Prop. \ref{prop:nsquare} severely limits the space of suitable $\lambda'$ within this chosen maximal torus $ T$ to a finite and discrete set of possibilities.

\subsection{Alignment, weight spaces and the permanent vs. determinant case} \label{subsec:weight-spaces}

Let us consider the case when $X\cong \C^{r+s}, W\cong \C^r$  and $f \in Sym^n (X^*)$ and $g\in Sym^n (W^*)$ be special forms with stabilizers $GL(X)_f=K$ and $GL(W)_{f_0}=H_W$. Let $\phi :W\rightarrow X$ be an invertible linear map such that the pull back of $f$ equals $g$, or in other words, $g = f\circ \phi$.  Let $Y=\phi (W) $ and $Z$ be a suitable complement of $Y\subseteq X$. Let us also identify $W$ with $Y$ and therefore $H_W$ with $H_Y$. Then, we can construct a $\lambda (t)\subseteq GL(X)$ such that we have the weight space decomposition $X=X_0 \oplus X_1$ with $X_0 =Y$ and $X_1 =Z$. We then have:
\[ \lambda (t^{-1}) f= t^0 f_0 + \ldots + t^m f_m \]
with $f_0 $ as the leading term and $f_0 |_Y =g$. The group $GL(X)_{f_0}=H$ has the following form:
\[ H= \left\{ \left[ \begin{array}{cc} 
a & b \\
0 & d \end{array} \right] | a \in GL(Z), b \in Hom(Z,Y), d\in H_Y \right\} \]
Let us suppose that there is an alignment between $f_0$ and $f$, i.e., a semisimple element $\mf{k} \in {\cal H} \cap {\cal K}^u$ as in Prop. \ref{prop:stabgeneral}.

\begin{prop} \label{prop:eigenspaces}
Suppose that $\mf{k}$ above has rational eigenvalues, then there is a $\phi' :W\rightarrow X$, a non-trivial 1-PS $\mu_X \subseteq K$ and a 1-PS $\mu_W \subseteq H_W$ such that:
\begin{enumerate} 
\item $\mu_X \circ \phi' = \phi' \circ \mu_W$. 
\item $f \circ \phi' =g$. 
\item If $X=\oplus_i X_i$ and $W=\oplus W_i$ is the weight space decomposition of $X$ and $W$ under $\mu_X$ and $\mu_W$, respectively, then $W_i = (\phi')^{-1}  (X_i )$, and thus $\phi' (W_i )\subseteq X_i$.
\end{enumerate}
\end{prop}

\noindent
{\bf Proof}: Let $\mf{k}\in  H \cap {\cal K}^u$ be a semisimple alignment and let $ \mu_X (a) =a^{\mf{k}}\subseteq GL(X)$. Then, since $\mf{k}$ commutes with $\lambda $, the 1-PS $\mu_X 
$ must preserve both $Y$ and $Z$. Let $\mu_Y $ be the restriction of $\mu$ to $Y$ and define $\mu_W = \phi^{-1} \circ \mu_Y \circ \phi$. Finally, let $\phi' = u \circ \phi$. It is easy to check (1) and (2). (3) follows from (1). $\Box $  

Let us apply this to the case where $X$ is a vector space of dimension $n^2$ with coordinate functions ${\cal X}=(X_{ij})_{i,j=1,\ldots ,m}$. 
Let $V=Sym^n (X^*)$ with $f(X)=det_n (X)$. Let $W$ be a space of  dimension $m^2+1$ with coordinate functions ${\cal W}=(W_{ij})_{i,j=1,\ldots,m} \cup W_{nn}$. Let $g_{m,n}=W_{nn}^{n-m} perm_m (W)$, the padded permanent.  

Suppose we have a $\phi :W \rightarrow X$ such that $f\circ \phi=g$ and the corresponding $\lambda $ and the partition $X=Y\oplus Z$ with $Y=\phi(W)$, such that: 
\[ \lambda (t) f = t^0 f_0 +t^1 + f_1 +\ldots + t^m f_m \]
such that $f_0 \circ \phi =g_{m,n}$. Prop. \ref{prop:eigenspaces} allows us to connect the weight spaces of stabilizer elements of the padded permanent with that of the determinant.

Towards this, we define: 
\begin{defn}
Let $A=(\{ 1,\dots , m\} \times \{1,\ldots , m\} )\cup (\{n\} \times \{n \})$ and $B=\{ 1,\ldots , n\} \times \{ 1,\ldots , n \}$ be sets of array indices. For a subset $R\subseteq A$, let $W_R =\{ w\in W | W_{i,j}(w)=0 \mbox{ for all } (i,j)\not \in R \}$. Thus, $W_R$ is the subspace of all vectors whose support is in the set $R$. Similarly, for $S\subseteq B$, we define $X_S$. A rectangular partition ${\cal R}=\{ R_1 ,\ldots , R_r \}$ of $A$ is where each $R$ is of the form $I_i \times J_j$, where $(I_i)$ and $(J_j )$ are two partitions of the row set and, respectively, the column set of $A$. Each rectangular partition ${\cal R}$ gives us a decomposition of $W=\oplus_{R\in {\cal R}} W_R$.  Similarly, we define a rectangular partition ${\cal S}=\{ S_1 ,\ldots , S_s \}$ of $B$ and the partition $X=\oplus_{S\in {\cal S}} X_S$.
\end{defn}

\begin{prop} \label{prop:rectangular}
Let $\lambda (t)$ be as above and suppose that there is a rational alignment between $f=det_n$ and $f_0 =Y_{nn}^{n-m} perm_n (Y)$. Then there is:
\begin{enumerate} 

\item a map $\phi'' :W \rightarrow X$ such that $g_{m,n}= det_n \circ \phi ''$
\item a rectangular partition ${\cal R}=\{ R_1 ,\ldots , R_r\}$ of $A$, and a rectangular partition ${\cal S}=\{ S_1 ,\ldots ,S_s \}$ of $B$, 
\item a correspondence $\Phi \subseteq {\cal R}\times {\cal S}$, such that $\phi'' (W_{R_i})\subseteq \oplus_{S_j \in \Phi (R_i)} X_{S_j}$. 
\end{enumerate}
\end{prop}

\noindent 
{\bf Proof}:  By Prop. \ref{prop:eigenspaces}, we have the 1-PS $\mu_X$ and $\mu_Y$ and a $\phi ' : W\rightarrow X$ such that $\mu_W \subseteq H_W$ and $\mu_X \subseteq K$. However, $\mu_X$ may not lie in the image of the standard torus $D_n \times D_n$ within $K_n$, the stabilizer of $det_n$. However, we may find a $k\in K_n$ such that $\mu_X^k $ does indeed lie within this torus. Define $\phi'' = k\circ \phi '$. The 1-PS $\mu_X^k$ gives us the rectangular partition ${\cal S}$ of $B$ and the weight space decomposition of $X$.

We also have the decomposition of $W=\oplus_{i} W_i$ by the weights of $\mu_W$. Now the connected part of $H_W$, the stabilizer of the padded permanent is a sub-torus of $(\C^*)^m\times (\C^*)^m \times\C^*$, where the action of $(\overline{\alpha})\times (\overline{\beta})\times \gamma $ is given by $W_{ij} \rightarrow \alpha_i \beta_j W_{ij}$ and $W_{nn}\rightarrow \gamma W_{nn}$. 
Thus the 1-PS $\mu_W \subseteq H_W$ indeed gives us a rectangular partition $A$ and a weight space decomposition ${\cal R}$ of $A$. For a given weight, say $d$, we define $\Phi_d = \{ (R , S )\in {\cal R}\times {\cal S}| wt_{\mu_{W}} (W_R )=wt_{\mu_X} (X_S)=d\}$ and $\Phi =\cup_d \Phi_d $. The condition that $\phi'' W_d \subseteq X_d$ then implies (3). $\Box $ 

\begin{remark}
The above proposition shows that a rational alignment $\mf{k}$ leads to specific information about the function $\phi$ and the support within the matrix $X$ for every coordinate $w_{ij}\in W$.

In general,  Prop. \ref{prop:eigenspaces} also indicates the importance of weight-spaces for 1-PS within stabilizers of forms and the coupling achieved when there is alignment. This is true even when the leading term is not of degree $0$. For example, in the case of $det_3$, the form $Q_1$ is a leading term of degree $0$ and the alignment between the weight spaces is evident. This is seen even for $Q_2$, which is not a leading term of degree $0$. 

For both, the permanent as well as the determinant, these rectangular spaces are also linear subspaces within their respective hypersurfaces. Such subspaces are of interest as the following proposition illustrates.  
\end{remark}

\begin{prop} \label{prop:lowerbound}
Let $H_m \subseteq \C^{m^2}$ be the hypersurface of $perm_m$. Suppose that there is a function $k(m)$ and a sequence of points $(x_m )$ for every $m$ such that $x_m \in H_m$, and the guarantee that dimension of any linear subspace $W\subseteq H_m $ containing $x_m$ is bounded by $k(m)$. Then, if $x_{nn}^{n-m}perm_m =\widehat{det_n}^{\lambda}$, for some $\lambda (t)\subseteq GL_{n^2} (\C)$, then $n>m^2-k(m)-1$.  
\end{prop}

\noindent 
{\bf Proof}: Let $\lambda $ and $Y\subseteq X=\C^{n^2}$ be as above so that a suitable conjugate $det_n^g$ of $det_n$, when restricted to $Y$ gives us the padded permanent $x_{nn}^{n-m} perm_m$. Now, $x_m \in Y$ be as above. 
Since $det^g_n (x_n)=x_{nn}^{n-m} perm_m (x_m )=0$, there is a linear subspace $Z\subseteq X$ containing $x_m$ of dimension $n^2-n$ such that $det_m^g (z)=0$ for all $z\in Z$. Let $W=Z\cap Y$. Then $x_{nn}^{n-m} perm_m (x_m )$ vanishes on $W$. 
By the hypothesis, we have $dim(W)\leq k(m)$. But we also have $dim(W)\geq (n^2 -n)+(m^2+1)-n^2 $. Combining the two, we get:
\[ m^2 -1 -n \leq k(m)\]
Thus $n^\geq m^2-k(m)-1$. $\Box$

\begin{remark}

The presence of alignment and its use in proving lower bounds was explored in \cite{LANDSBERG2017146}. They propose a stronger form where, for every (or a large fraction) of elements $h\in H_Y$, we have an element $g\in K$ which preserves $Y$ and matches $h$. They obtain exponential lower bounds for $m$ in terms of $n$ for such $\phi$. The rectangular partitions are seen in their implementation of $perm_3$ via $det_7$, (with $x_{77}=1$) is as given in ~\cite{grenet2011upper} and  is reproduced below:
\[ 
\left[ \begin{array}{ccccccc}
0 & 0 & 0 & 0 & x_{31} & -x_{32} & x_{33} \\
x_{11} & 1 & 0 & 0 & 0 & 0 & 0 \\
x_{12} & 0 & 1 & 0 & 0 & 0 & 0 \\
x_{13} & 0 & 0 & 1 & 0 & 0 & 0 \\
0 & -x_{22} & x_{21} & 0 & 1 & 0 & 0 \\
0 & -x_{23} & 0 & x_{21}& 0 & 1 & 0 \\
0 & 0 & -x_{23} & x_{22}& 0 & 0 & 1 \\\end{array} \right] 
\]
The partitions are, of course, $I=\{ 1\} \{ 2\} \{ 3\}\{ 7 \}$ and $J=\{ 1,2,3\}\{ 7\}$ for the permanent and $I=J=\{ 1\} \{ 2,3,4\} \{ 5,6,7\}$ for the determinant. 

The lower bound in Prop. \ref{prop:lowerbound} has already been shown by \cite{mignon2004quadratic} using the Hessian of a generic point on the hypersurface of the determinant, and a special point $x_n \in H_n$ as above. Our result requires us to compute $k(n)$. It has been mentioned here for its connection with weight spaces of stabilizer elements.  
\end{remark}

Finally, what about the case when there is no alignment?

\begin{defn}
Let $Z_0 = \{ \widehat{det_n^g}^{\lambda} | g \in GL(X)\}$ be the collection of leading terms obtained by applying the special 1-PS $\lambda $  to all elements of the $GL(X)$-orbit $O(det_n )$, and let $\overline{Z_0 }$ be its closure. We call these the {\bf co-limits} of $f_0 = g_{m,n}(Y)$.   
\end{defn}

We then have:
\begin{prop}
For any semisimple element $s\in K_n$, there is a $\overline{u}\in \overline{U(\lambda)}$, the unipotent radical of $H$, and a $u\in U(\lambda )$ such that $s^{u\overline{u}}$ stabilizes $f_0 '=\widehat{det_n^{\overline{u}}} \in Z_0$.  Moreover, there is an irreducible component $Z^i$ of $\overline{Z_0}$ containing both $f_0$and $f_0 '$.  
\end{prop}
The structure of $Z_0$ is discussed in Section \ref{subsec:CoLimit}.
We also conjecture, \ref{conj:stabilityalign}, that the $GL(Y)$-stability of $f_0$, the padded permanent, and the uniqueness of the form within $Sym^n (Y^*)$ for its stabilizer, point to a necessary alignment between $f_0$ and $det_n$. 

\subsection{A computation}
By Prop. \ref{prop:unistablegeneral}, semisimple elements in ${\cal P}(\lambda) \cap {\cal K}$ lead us to the presence of their conjugates within ${\cal H}$. This offers a significant insight into the alignment between ${\cal K}$ and ${\cal H}$. But what if no such elements exist? If ${\cal T}$ is a maximal torus of ${\cal K}$, what if ${\cal T}\cap {\cal P}(\lambda)=0$? The degree conditions of lemma \ref{lemma:dicho} lead to us to define the integers: 
\[ \delta (\mf{g},\mf{g'})=deg([\mf{g},\mf{g}'])-deg(\mf{g})-deg(\mf{g}')\] 
For any $\mf{g},\mf{g}'\in {\cal G}$, we have $\delta (\mf{g},\mf{g}')\geq 0$ and the dichotomy that either (i) $\delta (\mf{g},\mf{g}')= 0$, or (ii) $[\hat{\mf{g}}, \hat{\mf{g}'}]=0$. These lead to significant combinatorial constraints which we now illustrate. 

Assume that there is a subalgebra ${\cal L}\subseteq {\cal K}$ where ${\cal L}\cong sl_r$ for some $r>0$. This happens for example when $y$ is the determinant polynomial. Let ${\cal T}_r$ be an identified maximal torus and suppose that $deg(\mf{t})<0$ for all $\mf{t}\in {\cal T}_r$. We can then find a basis ${\cal C}$ such that $\hat{\cal C}$ generates $\hat{\cal T}_r \subseteq {\cal G}_{-}$. Let ${\cal C}=\{ K_i | i=1,\ldots ,r-1\}$ be such a basis. Note that such bases exist in an open set of ${\cal T}_r^{r-1}$.

In terms of this basis we then have ${\cal X}=\{ X_{ij} | 1\leq i\neq j \leq r \}\subseteq {\cal L}$, the collection of root vectors. Together ${\cal C}\cup {\cal X}$ form a basis for ${\cal L}$.  In terms of this basis, there are the standard Lie bracket relations, some of which are presented below:
\begin{equation} \label{eqn:lieb}
\begin{array}{rclr}
[K_i , K_j ]&=& 0 & (a)\\[0pt]
[X_{ij},X_{kl}]&=& 0 \mbox{ when $j\neq k$} & (b) \\[0pt]
 [ K_i ,X_{jk}] &=& c' X_{jk} & (c) \\[0pt]
[X_{ij}, X_{jk}]&=& cX_{ik} \mbox{ when $i\neq k$}& (d) \\[0pt]
[X_{ij}, X_{ji}]&=& K_{ij} =\sum_{k=i}^{j-1} K_k \mbox{ when $i<j$}& (e) \\
\end{array}
\end{equation}
for some $c\in \R, c\neq 0$ and $c' \in \R$, possibly zero. 

For all $i,j$, let  $d_{ij}=deg(X_{ij})=deg(\hat{X}_{ij})$ and $k_{ij} =deg(K_{ij})=deg(\hat{K}_{ij} )$ and note that $k_{ij} <0$. 

Let us now analyse $\hat{\cal L}_{\cal C}\subseteq \hat{\cal L}$, the subalgebra generated by the leading terms $\hat{\cal X} \cup \hat{\cal C}$ of the chosen basis. Note that that while $\hat{\cal T}_r \subseteq \hat{\cal L}_{\cal C}$ irrespective of the basis ${\cal C}$, the algebra $\hat{\cal L}_{\cal C}$ is determined by the choice of ${\cal C}$ and {\em need not equal} $\hat{\cal L}$.  

Conditions $(a)$ and $(b)$ give us $[\hat{K}_i ,\hat{K}_j]=0$ and $[\hat{X}_{ij} , \hat{X}_{kl}]=0$ when $j\neq k$. Looking at (c), we have $[K_i ,X_{jk}]=c'X_{jk}$ and comparing degrees, we see that the total degrees on the left and right do not match, i.e., $d_{jk}+k_i \neq d_{jk}$ and hence $[\hat{K}_i , \hat{X}_{jk}]=0$. In other words $[\hat{\cal T}_r, \hat{\cal L}_{\cal C}]=0$ and $\hat{\cal T}_r$ is in the center of $\hat{\cal L}_{\cal C}$. 

For any triple $r,s,t$, with $r\neq s$ and $s\neq t$, let $\delta_{rst}=d_{rt}-d_{rs}-d_{st}$ corresponding to the equation $[X_{rs},X_{st}]=cX_{st}$ (with $c\neq 0$). We then have $\delta_{rst}\geq 0$. 
Now consider a triple $i,j,k$ of distinct numbers, we have:
\[ \begin{array}{rcl} 
[ X_{ij}, X_{ji}]&=& K_{ik} \\[0pt]
[X_{ij}, X_{jk}] &=& X_{ik} 
\end{array}
\]
This gives is the degree conditions:
\[ \begin{array}{rcl} 
d_{ij}+d_{ji}+\delta_{iji}&=&k_{ik}\\
d_{ij}+d_{jk}+\delta_{ijk} &=& d_{ik} 
\end{array}
\]
Eliminating $d_{ij}$ and rearranging, we get:
\[  \begin{array}{rcl}
d_{ji}+d_{ik}-d_{jk}+\delta_{iji}&=& k_{ik}+\delta_{ijk}
\end{array}
\]
In other words:
\[ \delta_{jik}+\delta_{iji}=k_{ik}+\delta_{ijk}
\]
Since $k_{ik}<0$, this forces the condition $\delta_{ijk}>0$, or in other words $[\hat{X}_{ij}, \hat{X}_{jk}]=0$! What about $[\hat{X}_{ij} ,\hat{X}_{ji}]$? Since $[X_{ij}, X_{ji}]=K_{ij}$, we either have $[\hat{X}_{ij}, \hat{X}_{ji}]=0$ or $[\hat{X}_{ij}, \hat{X}_{ji}]=\hat{K}_{ij}$. 
Thus for all $i,j,k,l$, we have $[\hat{X}_{ij}, \hat{X}_{kl}]\in \hat\C \cdot {\cal C}$. Thus $\hat{\cal L}_{\cal C}/\hat{\cal T}_r$ is abelian. So $\hat{\cal L}_{\cal C}$ is an abelian extension of $\hat{\cal T}_r$. 

\section{The normal cone}\label{sec:normalcone}

In this section, we examine the role of $y_e$, the tangent of approach to the point $z$ and $O(z)$, and define suitable $G$-varieties $W\supseteq \overline{O(z)}$ which allow approaching $z$ along the tangent $y_e$ while staying within $W$.  This analysis considers a more general 1-parameter family (1-PF) $\gamma (t)\subseteq G$  for taking limits. This allows us to handle elements $z\in \overline{O(y)}$ which do not arise as leading terms of 1-PS. 

\subsection{Generalities}
\begin{defn}
$\C [[t]]$ will denote the ring of formal power series, and $\C ((t))$ the ring of Laurent series, the quotient field of $\C [[t]]$. A 1-parameter family (or simply 1-PF) $\gamma$ is a family of group elements $\gamma(t)=(g_{ij}(t))$, where each $g_{ij} (t)\in \C ((t))$.     
\end{defn}

\begin{remark}
By Theorem 1.4 \cite{kempf1978instability}, 1-PF above are adequate to detect closure in the Zariski topology. In other words, if $z\in \overline{O(y)}$, the Zariski closure, then there is a 1-PF $\gamma (t)$, where $\gamma(t)$ is a matrix in $G$ with power series entries, such that \[ y(t)=\gamma(t) y =y_0 +\sum_{i\geq 1} y_i t^i \]
with $z=y_0$.
\end{remark}

\begin{lemma}    
Let $\gamma $ be a 1-PF and let $v\in V$ and $v\neq 0$. Let $v(t)=\rho (\gamma(t))(v)$, then there is an $a\in \Z$ such that:
    \[ v(t)= \sum_{i\geq a} v_i t^i\]
    with $v_i \in V$ for all $i$ and $v_a \neq 0$. Similarly, for $\mf{g}\in {\cal G}, \mf{g}\neq 0$ and $\mf{g}(t)= adj(\gamma(t))(\mf{g})$, then there is a $b\in \Z$ such that:
\[ \mf{g} (t)=\sum_{i\geq b} \mf{g}_i t^i\]
with $\mf{g}_i \in {\cal G}$ for all $i$ and $\mf{g}_b \neq 0$. 
\end{lemma}

\begin{remark}
    Note that, unlike the case of a 1-PS, we may not have a decomposition of the ambient space $V$ as well as the Lie algebra ${\cal G}$. 
\end{remark}
\noindent 
{\bf Proof}: Suppose that $dim(V)=s$. From the rationality of the representation $\rho$, and for a chosen basis ${\cal B}$ of $V$, $\rho(\gamma(t))({\cal B})=A(t){\cal B}$, where ${\cal B}$ is $1\times s$ and $A(t)$ is an $s\times s$-matrix with entries in $\C ((t))$. If $v=c{\cal B}$, then $\gamma (t)\cdot v=v(t)=cA(t){\cal B}$. Thus, there is a row-vector $c(t)$ such that $v(t)=c(t){\cal B}$. Unpacking $c(t)$ by degrees gives us $v(t)=\sum_{i\geq a} (c_i {\cal B})t^i$, where $c_i \in \C^{1\times s}$ are row vectors. Using these as the coefficients for ${\cal B}$ proves that $v(t)=\sum_{i\geq a} v_i t^i $.

Coming to the second statement, let ${\cal B}=\{ \mf{g}_1 , \ldots , \mf{g}_r \}$ now be a basis for ${\cal G}$. 
Each matrix $\mf{g}_i \in {\cal G} \subseteq gl_n $ may be expressed as a column vector $\mf{c}_i \in \C^{n^2}$. By the same calculation as above we have $\gamma (t)\cdot \mf{g}=\mf{g}(t)$ is a matrix in $gl_n$ with entries in the Laurent series $\C ((t))$, and therefore as a column vector $\mf{c}(t)$ too. 
Let $A(t)$ be the $n^2 \times (r+1)$-matrix $[\mf{c}_1 ,\ldots , \mf{c}_r , \mf{c}(t)] $. 
Clearly, since $\gamma(t)\subseteq G$, we have $\gamma (t)\cdot \mf{g}\in {\cal G}$, and the rank of this matrix (with entries in $\C ((t))$) is $r$.  Therefore $\mf{g}(t)$ is expressible as $\sum_{i=1}^r a_i (t) \mf{g}_i $, for some elements $a_i (t)\in \C ((t))$. Collecting terms of the same degree gives us the result. $\Box $

\begin{defn}

For an element $\mf{k}\in {\cal K}$,  we define $\mf{k}(t)$ as the element $\gamma(t)\cdot \mf{k} \in {\cal G}\otimes \C ((t))$ and $\hat{\mf{k}}$ as the leading term of $\mf{k}(t)$. The space $\hat{\cal K}$ will denote the $\C$-space formed by all elements $\{ \hat{\mf{k}} |\mf{k}\in {\cal K} \}$. The space ${\cal K}(t)$ will denote the $\C ((t))$-space formed by the elements $\{ \mf{k}(t)| \mf{k}\in {\cal K} \}$. In other words, ${\cal K}(t) = {\cal K} \otimes C((t))$.
\end{defn}

\noindent 
Note that ${\cal K}(t)$ is a Lie algebra over $\C ((t))$. We then have the following: 
\begin{prop}
In the above notation, we have:
\begin{enumerate}
    \item There is a basis $\{ \mf{k}_i (t) \}_{i=1}^r \subseteq {\cal G}\otimes \C [[t]]$ of ${\cal K}(t)$ such that (i) $r=dim({\cal K})$, (ii) if ${\cal K}(0)$ is the space formed by the elements $\{ \mf{k}_i (0)\}_{i=1}^r$, then ${\cal K}(0)$ is a subalgebra of ${\cal H}_{\overline{y_e}}\subseteq {\cal H}$ of the same dimension. 
    \item $\hat{\cal K}\subseteq {\cal K}(0)$. 
\end{enumerate}
\end{prop}

\noindent 
The proof is a careful reworking of the proof of lemma \ref{lemma:samedim}, for details see \cite[Theorem 3.13] {liemethodsGCT}. We omit the proof.  Note that $\hat{\cal K}\subsetneq {\cal K}(0)$ is eminently possible.

We continue with the notation $\gamma (t) \subseteq G$ such that:
\[ y(t)=\gamma (t)y= z+y_e t^e + \mbox{ higher degree terms} \]
where $y_0 =z$ and $y_e$ is the tangent of approach. Note that for the 1-PS $\lambda $ of earlier sections, we do have the 1-PS $\lambda'$ which is of the above format. 

Recall that $I_y $ (resp. $I_z$) are ideals in $\C [V]$ for the varieties $\overline{O(y)}$ (resp. $\overline{O(z)}$). The rings $A_y$ (resp. $A_z $) are the corresponding coordinate rings, i.e., $\C [V]/I_y $ (resp. $\C[V]/I_z$). Note that since $\overline{O(y)}$ and $\overline{O(z)}$ are cones, the ideals $I_y$ and $I_z$ are homogeneous. 

We use this limiting 1-PF $\gamma$ to construct a suitable set of derivations ${\cal D}$ on the ideal $I_z$. We then use this to define a $G$-invariant ideal $\overline{J}$ in the associated graded ring of ${\mathbb C}[V]$ with respect to $I_z$. 

We begin with basic results on graded rings. 

\begin{defn}
 Let $R =\oplus_{i\geq 0}  R_i$ be the associated graded ring for the ideal $I_z \subseteq \C [V]$. In other words, $R_i =I_z^i/ I_z^{i+1}$ (with $I_z^0 =\C [V]$). For any homogeneous ideal $I$, 
 set $I_i = I^i_z \cap I$ and $\overline{I}_i= (I^{i+1}_z +I_i)/I^{i+1}_z$. Let $\overline{I}=  \oplus_{i\geq 0} \overline{I}_i$ be the filtration of the ideal $I$. 
\end{defn}

Note that $(I_i)$ is an $I_z$-stable filtration of the ideal $I$ and that $\overline{I}$ is an ideal within $R$. 
The ring $R$ is the ring of functions on the ``normal cone'' to the variety $\overline{O(z)}$. We have
\begin{prop} \label{prop:asso}
Using the above notation:
\begin{enumerate}
    \item The ring $\C [V]$ is isomorphic to $R$ as $G$-modules. For any $G$-invariant ideal $I\subseteq \C [V]$ $I$ and $\overline{I}$ are isomorphic as $G$-modules and so are $\C[V]/I$ and $R/\overline{I}$.  
    \item The ring $A_y$ is isomorphic to $R_y = R /\overline{I_y} =\sum_{i\geq 0} R_i/\overline{(I_y)}_i $ as $G$-modules. Moreover, $(R_y)_0 =A_z$. 
    \item We have the exact sequence of ideals and rings (as well as $G$-modules):
    \[ 0 \rightarrow \overline{I_z}/\overline{I_y} \rightarrow R/\overline{I_y} \rightarrow R/\overline{I_z} \rightarrow 0 \]
    In this sequence $(R/\overline{I_z})_i =0$ for all $i>0$, and thus $ (\overline{I_z}/\overline{I_y})_i \cong (R/\overline{I_y})_i$, for $i>0$. 
\end{enumerate}
\end{prop} 

\noindent
The proof is clear and is omitted. \\

\subsection{The tangent ideal $\overline{J}$}
The main result of this subsection is an extension of the Peter-Weyl condition and is given below. It is based on the construction of the ideal $\overline{J} \subseteq R$.  

\begin{prop} \label{prop:mainrep}
Let $z$ be the leading term of $y$ under the action of the 1-PF $\gamma $ and let $y_e$ be the tangent of approach. Let ${\cal H}_{\overline{y_e}} \subseteq {\cal H}$ be the stabilizer of $\overline{y_e} \in V/T_z O(z)$. Then for all $i\geq 1$, then the dual $(R_i /\overline{(I_y)}_i)^*$ is non-zero and has an ${\cal H}_{\overline{y_e}}$-fixed vector. 
\end{prop}

\begin{remark} 
For any subgroup $L\subseteq G$, let $Rep_G (L)$ be those $G$-modules $W$ such that $W^*$ has an $L$-fixed vector. 
The importance of this is in the Peter-Weyl theorem on closed $G$-orbits $O(v)$ with reductive stabilizer $L=G_v$. In this case, we have the isomorphism of $G$-modules $\C[O(v)]\cong \C[G]^L \cong  \oplus_{W \in Rep_G (L)} n_W W$, where $n_W$ is the dimension of the space of $L$-fixed vectors in $W^*$.  In the above notation, the kernel $I_z/I_y$  of the surjection $A_y \rightarrow A_z$ contains $G$-modules in the set $Rep_G ({\cal H}_{\overline{y_e}})$. The other modules in $A_y$ come from $A_z$ and belong to $Rep_G ({\cal H})$
\end{remark}

The proof of the above proposition needs the construction of the ideal $\overline{J}$, which we now proceed to do.

\begin{defn}
Let $w\in V$ be a point, $T_w V$ be the tangent space at $w$ and $v\in T_w V$ be a tangent vector at $w$. For an indeterminate $\epsilon$, we consider the substition and the expansion:

\[  f(w+\epsilon v )=f_0 + f_1\epsilon + \ldots \]
where $f_i \in \C$. We define $D_{w,v} : \C [V] \rightarrow \C $ as $D_{w,v} (f)=f_1$. 
\end{defn}

The functional $D_{w,v}$ is called a derivation at $w$. It is easy to check that for any $f,f'\in \C [V]$, we have $D_{w,v} (ff')= f(w)D_{w,v} (f')+ f'(w) D_{w,v} (f)$. 

For the $\gamma (t) \subseteq G$ and $g\in G$, we define $\gamma^g (t)$ as the family $g\gamma(t)\subseteq G$.  Let $y(t)=\gamma (t)y$ and $y^g (t)=\gamma^g (t)y$. Then we have:
\[ y^g (t)= gz+(gy_e)t^e + \mbox{ higher degree terms } \]
Note that both $y(t)$ and $y^g (t)$ are elements of $V\otimes \C [[t]]$, i.e., power series in $t$ with coefficients in $V$. 

\begin{lemma} \label{lemma:dequiv}
Given any polynomial $f\in \C [V]$ and for any $g\in G$, we make the substitution as below:

\[ f(y^g (t))=f_0 +f_1 t^1 + \ldots \]
Then $f_0 =f(gz)$, $f_1= \cdots = f_{e-1}=0$ and $f_e =D_{gz,gy_e} (f)$.
\end{lemma}

\noindent
The proof is a simple computation. This makes the substitution of $y^g(t)$ within $f$ behave as a path which evaluates the derivative of $f$ at $gz$ along the direction of the tangent of approach $gy_e$.  

\begin{defn} \label{defn:Dk}
Let $w\in \overline{O(z)}$, $v \in T_w V$ be arbitrary. For a $k\geq 1$, and an element $f\in I_z^k$, upon substitution of $w+\epsilon v$ in $f$, we have:
\[ f(w+\epsilon v)= f_0 +f_1 \epsilon \ldots \]
Then it is clear that $f_0 = \cdots = f_{k-1}=0$. We define $D^k_{w,v} (f)$ as $f_k$, the coefficient of $\epsilon^k$. 
\end{defn}

$D^k_{w,v}$ is well defined and  $D^k_{w,v} (f)=0$ for $f\in I_z^{k+1}$, whence $D^k_{w,v}$ is effectively a map $D^k_{w,v} : I^k_{z}/I^{k+1}_z \rightarrow \C$. 

\begin{prop} \label{prop:mainderivation}
For the notation as above, let $f\in I_z^k$ and $g\in G$, then $f(y^g (t))$ is a power series in $t$ of the form: 
\[ f(y^g (t))=  \alpha t^{ke} + \mbox{ higher terms in $t$} \]
where $\alpha =D^k_{gz,gy_e} (f)$. 
\end{prop}

\noindent
The proof is straightforward. 

\begin{prop} \label{prop:basicDk} 
With the above notation, for any $w\in \overline{O(z)}, v\in T_w V$,  we have:
\begin{enumerate} 
    \item For $f, f'\in I^r_z$ and $f''\in I^s_z$, we have:
\[ \begin{array}{rcl}
D^k_{w,v} (f)&=& 0 \mbox{ for all $k<r$} \\
D^r_{w,v} (f+f')&=&D^r_{w,v} (f)+D^r_{w,v} (f')\\
    D^{r+s}_{w,v} (ff'')&=&D^r_{w,v} (f) D^s_{w,v} (f'') \end{array}\]
\item For any $g\in G$, and $f\in I^r_z$, let $f^{g}$ denote the function $f^g(x) = f(g^{-1} x)$. Then, $f^g \in I^r_z$ as well. Moreover, for any $g'\in G$, we have:
\[ D^r_{gz,gy_e}  f^{g'} =D^r_{(g')^{-1} gz, (g')^{-1}gy_e } (f) \]
\end{enumerate}
\end{prop}
    
\noindent 
{\bf Proof}: The proof of (1) is straightforward. 
For (2), note that $f^{g'}(y^g(t))=f((g')^{-1} y^g(t))=f(y^{(g')^{-1}g} (t))$. This proves the proposition. $\Box $

\begin{lemma} \label{lemma:gomod}
For $w\in O(z)$ and $v\in T_w V$ let $v' \in T_w O(z)$, the tangent space of the orbit $O(z)$ at $w$. Then for an $f\in I^r_z$ with $r>0$, we have $D^r_{w,v+v'}(f)=D^r_{w,v} (f)$.     
\end{lemma}

\noindent
{\bf Proof}: Let $f' \in I_z$. Since $v'\in T_w O(z)$, the tangent space of the orbit $O(z)$ and since $f'$ vanishes on $O(z)$, it is easy to check that $D_{w,v'} (f')=0$. Hence, for any $f\in I^r_Z$, we have $D^r_{ w,v'} (f)=0$ as well. This proves the claim. $\Box $

\begin{remark}\label{rem:yebar}
Proposition $\ref{prop:basicDk}$ tells us that since $D^k_{w,v}$ vanishes on $I_z^{k+1}$, it is effectively a functional on $I_z^k /I_z^{k+1}$. The above lemma tells us that the functional depends on the representative $\overline{v} \in \overline{N}=V/T_w O(z)$, and not on $v$ itself.  
\end{remark}

We now use the tangent of approach $y_e$ to define the modules $\overline{J}_k $. 

\begin{defn}
For the above data, let $\overline{J}_k \subseteq R_k$ be defined as follows:
\[ \overline{J}_k =\{ \overline{f} \in I^k _z /I^{k+1}_z | D^k_{gz,gy_e} (\overline{f})=0 \mbox{ for all $g\in G$}\} \]
\end{defn}

\begin{prop} \label{prop:Jlambda}
With the above definition of $\overline{J}_k$, we have:
\begin{enumerate} 
\item $\overline{J}_r  I^s_z \subseteq \overline{J}_{r+s}$ and 
\item $(\overline{I}_y)_k \subseteq \overline{J}_k$.
\item Let $\overline{J}_{z,y_e}$ or simply $\overline{J}=\oplus_{k\geq 1} \overline{J}_k$. Then $\overline{J}$ is a $G$-invariant ideal of $R$ and  $\overline{I_y} \subseteq \overline{J}\subseteq \overline{I_z}=\oplus_{i\geq 1} R_i $. 
\end{enumerate}
Henceforth, we call $\overline{J}=\overline{J}_{z,y_e }$ as the {\bf tangent ideal} for the data $(z,y_e )$. 
\end{prop}

\noindent 
{\bf Proof}: For (1), we note that for any $f\in \overline{J}_r $ and $f'\in I^s_z$, have 
\[ D^{r+s}_{gz,gy_e} (ff') = D^r_{gz,gy_e} (f) D^s_{gz,gy_e} (f')\]
Since the first term is zero, so is the product. For (2), let $f''\in I_y \cap I^k_z$. 
Note that, by Prop. \ref{prop:mainderivation}, $D^k_{gz,gy_e} (f'')$ is also obtained by evaluating $f''$ on the path $y^g (t)=\gamma^g (t)y$. Since $y^g (t)\in O(y)$ for all $t>0$, we have $f''(y(t))=0$ for all $t>0$. Passing to the limit, we must have $D^k_{gz,gy_e }(f'')=0$. This proves (2). For (3) above, (1) ensures that $(\overline{J}_k)$ is $I_z$-stable. 
Its $G$-invariance follows from the fact that if for some $f\in I_z^k$,   $D^k_{gz,gy_e} (f)$ vanishes for all $g$, then, by Prop. \ref{prop:basicDk} (2),  so do $D^k_{gz,gy_e} (f^{g'})$, for any $g'$-translate of $f$. 
This proves the first part of (3). The second part follows from (2). This completes the proof of the proposition. 
$\Box $ 

\begin{remark} 
Thus, the variety $\overline{W} \subseteq Spec(R)$  of $\overline{J}$ is an infinitesimal $G$-invariant thickening of $O(z)$ within the normal bundle, and in the direction $y_e$ and its translates. Assuming that $y_e \not \in T_z O(z)$, and in light of Remark~\ref{rem:yebar}, this thickening is in a direction transverse to the tangent space to the orbit at $z$. A natural question is whether there is an ideal $J'\subseteq \C[V]$ such that its filtered version $\overline{J'}$ equals $\overline{J}$? In that case, the variety $W'\subseteq V$ of $J'$ would be a model for $\overline{W}$ in the normal bundle. 
By Prop. \ref{prop:Jlambda}, this model would be an intermediate variety and lie between $\overline{O(y)}$ and $\overline{O(z)}$. The construction of $J'$ and a recipe for computing its dimension are addressed in the next subsection. 
\end{remark}

We present three examples which illustrate that with $y\stackrel{\gamma}{\rightarrow} z$, the ideal $\overline{J}$ and the existence of an intermediate variety depend crucially on $y_e$. Moreover, the construction allows for $\overline{O(z)}$ to be singular within $\overline{O(y)}$. Also see Ex. \ref{ex:popov} which will be discussed later.

\begin{ex}
Let $V=\C^2$, and $z=[0,0]^T$ and $y=[1,1]^T$.  Let $G$ be the parabolic group given below. $V$ consists of $3$ orbits, viz., $O(z)$ of the sole point $z$, $O(y)$ of $\C^2$ minus the $X$-axis, and $O([1,0]^T)$ of the $X$-axis minus $z$. The stabilizer $H$ of $z$ is $G$. 

Consider the two families $\gamma_1 $ and $\gamma_2 $ below and the paths $\gamma^1 (t)y$ and $\gamma^2(t)y$ which take $y$ to $z$ but with tangents $y_1 =[1,1]^T$ and $y_2 =[1,0]^T$. 
\[ G=\left\{ \left[\begin{array}{cc}
a & b \\
0 & c \end{array} \right] | ac\neq 0 \right\} \; \: 
\gamma^1 (t) = \left[\begin{array}{cc}
t & 0 \\
0 & t \end{array} \right]
\: \: 
\gamma^2 (t)=\left[\begin{array}{cc}
t & 0 \\
0 & t^2  \end{array} \right]
\]
Let $\overline{J}_1$ and $\overline{J}_2$ be the tangent ideals for the data $(z,y_1 )$ and $(z,y_2 )$ respectively. Then $\overline{J}_1$ consists of all forms $f$ of degree $k$ whose $k$-th derivatives vanish in the directions $[1,1]^T$ and its $H$-orbit. 
Thus $(\overline{J}_1)_k=(0)$ for all $k$. Let $J_1 =(0)\subseteq \C[V]$, then $J_1$ is an ideal and $\overline{(J_1)}=\overline{J}_1$ and thus $W_1=\C^2$ is the required model. This does lie between $\overline{O(z)}$ and $\overline{O(y)}$ but is not strictly intermediate. 

On the other hand, since the orbit of $y_2=[1,0]$ is its own multiples, $(\overline{J}_2)_k$ is the ideal generated by the forms $\{ y^i x^{k-i}| i>0\}$. Again, there is the ideal $J_2 =(y)\subseteq \C[V]$ such that $\overline{(J_2)}=\overline{J}_2$. The variety $W_2 =\overline{O([1,0]^T)}$ is the required model.  Note that in this case $\overline{O(z)} \subsetneq W_2 \subsetneq \overline{O(y)}$ and $W_2$ is strictly intermediate.
\end{ex}

\begin{ex}
Let us consider the group $G=S^1 \times \R^*$ acting on $\R^3$ given by the matrix $g(\theta , r)$ as:
\[ g(\theta ,r)=\left[ \begin{array}{rcl}
r \cos \theta & r \sin \theta & 0 \\
-r \sin \theta & r \cos \theta & 0 \\
0 & 0 & r \end{array} \right] \]
Consider the point $y=[1,0,1]^T $ and its orbit $O(y)= \{ (r \cos \theta , r \sin \theta ,r)| \theta \in [0,2 \pi), r \neq 0 \}$. 
The point $z=[0,0,0]^T $ is in $\overline{O(y))}$ via $\lambda (t)=diag(t,t,t)$, with $y_e =[1,0,1]$. Let us compute the tangent ideal $\overline{J}$ for this data. 

Note that $O(z)=z, G_z=G$ and $I_z =(x,y,z)$. Let $f_1 =ax+by+cz \in I_z$, for some real $a,b,c$, $f_2 =x^2 \in I_z^2 $ and $f_3 =x^2+y^2-z^2 \in I_z^2$. Note that $\overline{O(z)}=z$ and $G_z =H$. We have:
\[ \begin{array}{rcl}
D^1_{z,y_e } (f_1 )&=& a+c \\
D^2_{z,y_e} (f_2 )&=& 1  \\
D^2_{z,y_e}  (f_3) &=& 0  \\
\end{array} \]
Applying a general group element $g(\theta ,r)$ and evaluating $D^k_{gz,gy_e }(f_i)$ gives us that $f_1 \not \in \overline{J}_1$ for any non-zero tuple $(a,b,c)$. Indeed $\overline{J}_1 =(0)$. On the other hand, $f_2 \not \in \overline{J}_2 $ but $ f_3 \in \overline{J}_2 $. In fact, for the ideal $J' =(x^2+y^2-z^2)$, we see that $\overline{J'}=\overline{J}$. This gives $\overline{O(y)}$ as $W'$ the intermediate variety. Note that $W'$ has a singularity at $z$. 
\end{ex}

We give another example of $\overline{J}$ from the classical representation of matrices under the adjoint action.

\begin{ex}
Let us consider $V=\C^{3 \times 3}$ with the coordinate functions $X=(X_{ij})$ in $(\C^{3\times 3})^*$. Consider the action of $GL_3 (\C)$ on $\C^{3\times 3}$ by conjugation. Thus, for a matrix $A$ and $g\in GL_3 (\C)$, the action of $g$ on $A$ is given by $A\rightarrow gAg^{-1}$. Consider the matrix $C$ and the family $\lambda (t)$ below:
\[ C=\left[ \begin{array}{ccc}
1& 0 & 1 \\
c_{21}& 1 & 0 \\
c_{31} & c_{32} & 1 \end{array} \right] 
\hspace*{0.8cm} \lambda (t) =\left[ \begin{array}{ccc}
1 & 0 & 0 \\
0 & t & 0 \\
0 & 0 & t^2 \end{array} \right] 
\]
The entries $c_{ij} \in \C$ are chosen such that the matrix $C$ has distinct eigenvalues. Let $N_1 $ , $N_2$ and $I $ be the matrices below:
\[ N_1 =\left[ \begin{array}{ccc}
0 & 0 & 1 \\
0 & 0 & 0 \\
0 & 0 & 0 \end{array} \right] 
\hspace*{0.4cm} N_2 =\left[ \begin{array}{ccc}
0 & 1 & 0 \\
0 & 0 & 1 \\
0 & 0 & 0 \end{array} \right] 
\hspace*{0.4cm} I =\left[ \begin{array}{ccc}
1 & 0 & 0 \\
0 & 1 & 0 \\
0 & 0 & 1 \end{array} \right] 
\]
We see that:
\[ \lambda(t)\cdot C= t^{-2} N_1 + I + \mbox{ terms with higher degree in $t$}\]
Thus, in the notation of this section, we have $y=C, z=N_1$ and $y_e=I$. 
The ideal of orbit closure $\overline{O(N_1 )}$ is given by the equation $X^2=0$. Moreover, $\overline{O (C)}$ is given by $O (C) \cup O (N_2 ) \cup O (N_1 ) \cup \bf{0}$, where $\bf{0}$ is the zero matrix. The action of the family $\gamma (t)=t^{2} \lambda (t)$ is given as:
\[ \gamma(t)y= N_1 +t^2 I + \mbox{ terms with higher degree in $t$}\]
Let us evaluate $D_{z,y_e} (X^2)$ and $D_{z,y_e} (X^3)$. Since $X^2 (N_1 )=X^3 (N_1)=0$ these are elements of $I_z$, and we have:
\[ \begin{array}{rcl}
D_{z,y_e} (X^2)
&=& 2N_1 \neq 0 \\
D_{z,y_2} (X^3) &=& N_1^2 I + N_1 I N_1 +I N_1^2 =0\\
\end{array}
\]
Thus $X^3 \in \overline{J}_1 $. This points to the variety
$W^1 =\overline{O (N_2)}$ as the possible intermediate variety. Note that $\overline{O (z)} \subsetneq W^1 \subsetneq \overline{O(y)}$. 
\end{ex}

We conclude this section with a proof of Proposition~\ref{prop:mainrep}.

\begin{prop} \label{prop:finalphi}
In the notation of Prof. \ref{prop:asso}, for 
each $i\geq 1$, $R_i /\overline{(I_y )}_i \neq 0$.  
\end{prop}

\noindent 
{\bf Proof}: Recall that $R_y = \oplus_{i\geq 0} R_i /\overline{I_y }_i $ is isomorphic to $A_y$ as a $G$-module. Since $I_y \subset I_z$, the $0$-th term, $(R_y)_0=\C[V]/I_z$ is precisely $A_z$. The kernel of 
$D^i_{z,y_e} : R_i /\overline{(I_y )}_i \rightarrow \C$ is precisely $\overline{J}_i$. Thus, it suffices to show that $R_i /\overline{J}_i$ is non-zero. Now, by the non-singularity of the point $z\in \overline{O(z)} \subseteq V$, there is an $f \in I_z$ such that $D_{z,y_e} (f)\neq 0$. Then $f^k \in R_i$ and $D^k_{z,y_e} (f^k)\neq 0$. 
$\Box $

\noindent 
{\bf Proof of Proposition \ref{prop:mainrep}}. In view of the above proposition, what remains to be shown is the existence of a non-zero functional on $R_i /\overline{(I_y)}_i$ which is ${\cal H}_{\overline{y_e }}$-invariant. This is of course $D^i_{z,y_e} : R_i /\overline{(I_y )}_i \rightarrow \C$. 
Indeed, for any $\mf{h}\in {\cal H}_{\overline{y_e }}$, we have $\mf{h} D^k_{z,y_e} =D^k_{z,\mf{h}y_e} $. Now since $\mf{h} y_e \in T_z O(z)$, for any $f\in I_z^k$, $D^k_{z,\mf{h}y_e} (f)=0$. This proves the invariance of $D^k_{z,y_e}$ and the proposition. $\Box $

\subsection{The structure of the tangent ideal}

This section analyses the tangent ideal $\overline{J}_{z,y_e}$ (or simply $\overline{J}$) and the question of its dimension.


Let $B$ be the affine variety $\overline{O(z)}\times V$ and $\C[B]$ its coordinate ring. Let $B'$ be the open subset $=\{ (w,v) \in B \: with \: w\in O(z)\}$ of $B$. For any $(w,v)\in B$, let us define the map $e(w,v):R\rightarrow \C$ as follows. For $\overline{f}=\sum_{i=0}^k \overline{f}_i \in R$, where $\overline{f}_i \in R_i$,we define:
\[ e(w,v)(\overline{f})= \sum_{i=0}^k D^i_{w,v} \overline{f}_i \]
Thus $e(w,v)$ treats $v$ as a member of the tangent space $T_w V$. We now list certain properties of $E=\{ e_{w,v} | w\in \overline{O(y)}, v\in T_w V \}$.

\begin{prop} \label{prop:mapphi}
For $B$ as above, we have:
\begin{enumerate}
    \item For any $(w,v)\in B'$, the map $e_{w,v}$ is a non-trivial algebra homomorphism. Hence, the kernel $M_{w,v} \subseteq R$ is a maximal ideal of $R$. Thus, we have a map $\phi :B' \rightarrow MaxSpec (R)$. 
    \item For any point $w\in O(z)$, if $v' \in T_w O(z)$ then $\phi (w,v)=\phi(w,v+v')$. Thus, the fiber of the map $\phi$ at any point in the $Im(\phi)$ is a linear space of dimension $dim(G)-dim(H)$.
    \item Let $B_J =\{ (gz,gy_e +v')| g\in G \mbox{ and } v' \in T_{gz}O(z) \}$. Then have:  
    \[ \overline{J} =\cap_{(w,v)\in B_J } M_{w,v} \]
\end{enumerate}
\end{prop}

\noindent 
{\bf Proof}: That $e(w,v)$ is an algebra homomorphism follows from Prop. \ref{prop:basicDk}. It is then clear that its kernel must be a maximal ideal of $R$. Hence (1) is clear. The second assertion follows from Lemma \ref{lemma:gomod}, which says that $D^k_{w,v+v'}=D^k_{w,v}$ when $v'\in T_w O(z)$. 
The fiber is clearly $dim(T_w O(z))$ which is $dim(G)-dim(H)$. 

Finally, coming to (3), It is clear that $\overline{J} \subseteq \cap_{(w,v)\in B_J} M_{w,v}$. In the other direction, note that $\lambda (t)\in H$ but $\lambda (t)y_e =t^{e-d} y_e$. Thus, not only is $(z,y_e)\in B_J$, but so is $(z,\alpha y_e )\in B_J$, for any $\alpha \in \C^*$. In general, for any $\alpha \in \C^*$, and $ (w,v)\in B_J$ we have $(w,\alpha v)\in B_J$ as well. Now:
\[ e(w,\alpha v) (\overline{f})=\sum_i \alpha^i e(v,w)\overline{f}_i \]
Whence $e(v,w)(\overline{f})=0$ implies that $e(v,w)(\overline{f}_i) =0$ as well. Thus for any  $\overline{f} \in \cap_{(w,v)\in B_j} M_{w,v}$, we have $\overline{f}_i \in \cap_{(w,v)\in B_J } M_{w,v}$ as well. But that is equivalent to the requirement that $\overline{f}_i \in \overline{J}_i$. Thus $\cap_{(w,v)\in B_J} M_{w,v} \subseteq \overline{J}$. This proves (3) and the proposition. $\Box $

\begin{prop}
    The dimension of $\overline{J}_{z,y_e}$ (i.e., $\overline{J}$) is $dim(G)-dim(H_{\overline{y_e}})$. 
\end{prop}

\noindent
{\bf Proof}: We have the map $\phi :B' \rightarrow MaxSpec(R)$ whose fiber at each point is of dimension $dim(G)-dim(H)$. By Prop. \ref{prop:mapphi} (3) above, $\phi^{-1} (MaxSpec(\overline{J}_{z,y_e}))=B_J$ whose dimension is $dim(G)-dim(H_{\overline{y_e}})+dim(G)-dim(H)$. This implies that the dimension of $\overline{J}_{z,y_e}$ must be $dim(G)-dim(H_{\overline{y_e}})$. $\Box$

We know in general that $\hat{\cal K}\subseteq {\cal H}_{\overline{y_e}}$. By the above theorem, if $dim({\cal K})< dim({\cal H}_{\overline{y_e}})$, then, in the normal cone $Spec(R)$, there is indeed an intermediate variety between $\overline{O(z)}$ and $\overline{O(y)}$. 

\begin{conjecture}
If $dim({\cal K})< dim({\cal H}_{\overline{y_e}})$, then there is a strictly intermediate variety $\overline{O(z)}\subsetneq W \subsetneq \overline{O(y)}$ of dimension $dim(G)-dim({\cal H}_{\overline{y_e}})$.     
\end{conjecture}

Next, we examine if there is an intermediate ideal $J'\subseteq \C [V]$ such that $\overline{J'}=\overline{J}$ and $I_z \supsetneq J'\supsetneq I_y$. Whether the ideal of the variety $W_{z,y_e}$ defined below works is not clear to us.

\begin{defn}
For any closed $G$-variety $W$ with $z\in W$, we say that $w_e$ is a tangent of approach at $z$ if there is a $w\in W$ and a 1-PF $\beta (t)\subseteq G$ such that:
\[ w(t)=\beta(t)w=z +w_e t^e + \mbox{ higher degree terms } \]
The collection of all tangents of approach is denoted by $\overline{T}_z W$ 
\end{defn}

\begin{defn}
Let ${\cal W}_{z,y_e}$ be the collection of all homogeneous $G$-varieties $W$ within $V$ which (i) contain $z$ (and therefore $\overline{O(z)}$) and (ii) for which $y_e \in \overline{T}_z W$. 
\end{defn}

\begin{lemma}
 Let $W_{z,y_e}= \cap_{W\in {\cal W}_{z,y_e}} W$. Then $W_{z,y_e}$ is an algebraic variety. Let 
 \[ I_{z,y_e}= \sum_{W\in {\cal W}_{z,y_e}} I_W\]
where $I_W \subseteq \C [V]$ is the ideal of $W$. Then  $\overline{I_{z,y_e}} \subseteq \overline{J}$. 
\end{lemma}

\noindent
{\bf Proof}: It is clear that $W_{z,y_e}$ is an algebraic variety and that $I_W $ is an ideal whose variety is $W$. Let us first show that for any $W$ as above $\overline{I_W}\subseteq \overline{J}$. Suppose that $f\in I_W \cap I_z^k$. By the definition of $W$, there is a 1-PF $\gamma(t)$ taking some $w \in W$ to $z$ with the tangent of approach being $y_e$. Plugging $w(t)=\gamma(t)w$ into the expression of $f$, we have 
\[ f(w(t)) = t^{ke} D^k_{z,y_e } (f) + \ldots \] 
Since $W$ is a $G$-variety, $w(t)$ lies entirely in $W$. 
Hence, $f \in I_W$ implies that $f(w(t))$ is identically zero, i.e., $D^k_{z,y_e} (f)=0$. The same applies to $D^k_{gz,gy_e }$ too. Hence $\overline{f}\in \overline{J}$. Thus $\overline{I_W}\subseteq \overline{J}$ and indeed $\overline{I_{z,y_e}}\subseteq \overline{J}$. 
$\Box $

\section{Intermediate strata and co-limits space} \label{sec:StrataCoLimit}

As in the previous section, we continue with the search for intermediate varieties. We go back to the 1-PS case with $\lambda (t)$ and its action as below: 
\[ \lambda (t)\cdot y = y_d t^d + y_e t^e +\ldots + y_D t^d \]
with $z=y_d$. 
In the first subsection, we use $T$, a maximal torus of $G$ containing $\lambda $, and elementary polyhedral theory to arrive at possible intermediate varieties. In the second subsection, we consider {\bf co-limits} of $z$, i.e., leading terms $\hat{y'}$ of degree $d$ for some $y' \in O(y)$. 

\subsection{Intermediate polyhedral strata} \label{subsection:character}
Let $H_0 = L(\lambda )\cap H$ and $T_H$ be a maximal torus in $H_0$. Since $\lambda (t)$ commutes with $T_H$, we may choose a maximal torus $T \subseteq G$ containing $\lambda$ as well as $T_H$. Let its rank be $r$ which equals the rank of $G$. We assume that $T\subseteq D_n \cong (\C^* )^n$, the group of diagonal matrices $diag(t_1 ,\ldots , t_n )$ in $GL(X)$. 

Let the general element of $T$ be $\overline{t}=(t_1 ,\ldots ,t_r ) \in (\C^* )^r$ and $\rho:T \rightarrow D_n$ be its realization within $GL(X)$. 
Through $\rho$, we have a finite set $\Xi (V)$ of characters, i.e., vectors $\chi =(\chi(i))$ in $\Z^r$,  and a weight space decomposition $V=\sum_{\chi \in \Xi (V)} V_{\chi}$ such that for any $v_{\chi} \in V_{\chi}$, we have $\rho (\overline{t})v_{\chi} =(\prod_{=1}^r t_i^{\chi(i)}) v_{\chi}$. The product is simply denoted as $\overline{t}^{\chi}$. For any $v\in V$ with $v=\sum_{\chi \in \Xi (V)} v_{\chi}$,  we define $\Xi (v)=\{ \chi | v_{\chi} \neq 0 \} $. 

Let ${\cal T}$ be the Lie algebra of $T$ and $\phi ({\cal T})  \subseteq {\cal G} \subseteq gl(X)$ be its realization. Let $(\mf{t}_1 , \ldots , \mf{t}_r )$ be a basis for ${\cal T}$ such that $t^{\mf{t}_i}=(1,1,\ldots, 1,t_i ,1\ldots ,1) \in (\C^*)^r \cong T$. 
Let ${\cal T}_{\Z}$ (resp. ${\cal T}_{\R}$) be the $\Z$-module  (resp. $\R$-module) generated by $\{ \mf{t}_1 ,\ldots , \mf{t}_r \}$. Henceforth, by ${\cal T}$, we will mean ${\cal T}_{\R}$. 

We fix an element ${\bf 1}\in {\cal T}_Z$ such that $\phi (t^{\bf 1}) \in Z$, the center of $GL(X)$.  
For a $\mf{t} \in {\cal T}_{\Z}$, there are $(\mf{t}(i)) \in \Z^r $, the coefficients of $(\mf{t}_1 ,\ldots , \mf{t}_r )$ such that $t^{\mf{t}}=(t ^{\mf{t}(1)},\ldots , t^{\mf{t}(r)}) \in T$. 
Moreover, we have $\chi (t^{\mf{t}})=t^{\langle \chi , \mf{t} \rangle}$, where $\langle ,\rangle$ is the usual inner product on $\R^r$. For any 1-PS $\mu (t)\subseteq T$, we have $\log (\mu)\in {\cal T}_{\Z}$ such that $t^{\log (\mu)}=\mu(t)$. Finally, for $\bf{1}\in {\cal T}_Z$ there is a $c\neq 0$ such that for all $\chi \in \Xi (V)$, $\langle \chi , {\bf 1}\rangle =c$

\begin{lemma} 
In the above notation, for any $v\in V$, with $v=\sum_{\chi \in \Xi (v)} v_{\chi} $ and $\mf{t}\in {\cal T}$, we have:
\[ {\mf{t}} v =\sum_{\chi \in \Xi (v)} \langle \chi , \mf{t} \rangle v_{\chi} \]
\end{lemma}

\begin{defn} For any $v \in V$, let ${\cal T}_v =\{ \mf{t}\in {\cal T} | \mf{t}v=0\}$ be the Lie algebra stabilizer of $v$ within ${\cal T}$. 
For our $y\in V$, let ${\cal T}^+(y)=\{ \mf{t} \in {\cal T} | \langle \chi , \mf{t} \rangle \geq 0 \mbox{ for all } \chi \in \Xi (y) \}$. For any $w\in V$ let ${\cal T}^+_w$ be ${\cal T}_w \cap {\cal T}^+ (y)$. 
\end{defn}

\begin{lemma}
For the above data, let
\[ {\cal T}_{\Z , v} = \{ \mf{t}  \in {\cal T}_{\Z} | \langle \chi , \mf{t}\rangle =0 \mbox{ for all } \chi \in \Xi (v) \} \] 
Then ${\cal T}_v $ equals ${\cal T}_{\Z ,v} \otimes_{\Z} \R$. 
\end{lemma}
\begin{defn}
A set $F\subseteq \Xi (y)$ is called a face if there is a $\mf{t}\in {\cal T}^+ (y)$ such that $F=\{ \chi \in \Xi (y)| \langle \chi , \mf{t} \rangle =0 \}$. This face $F$ is also denoted by $F(\mf{t})$. The dimension $dim(F)$ of a face is the dimension of the vector space $\R \cdot F$ within $\R^n$ spanned by the elements of $F$. A face $F$ of co-dimension $1 $ within $\Xi (y)$ is called a facet. The element $y_F$ is defined as $\sum_{\chi \in F} y_{\chi}$. 
\end{defn}

Note that $\Xi (y)$ is also a face and that $\Xi (y_F )=F$. We have the following simple lemma:

\begin{lemma}\label{lemma:stratbasic}
For the above data, we have:
\begin{enumerate} 
\item Let $F \subseteq \Xi (y)$ and $y_F$ be as above. Then the dimension of the face $F$  complements the dimension of the sub-torus ${\cal T}_{y_F}\subseteq {\cal T}$ which stabilizes it, i.e., $dim(F)=r-dim({\cal T}_{y_F} )$.
\item Let $\mf{t} \in {\cal T}^+ (y) $ and $\mu=t^{\mf{t}}$. Then the leading term $\hat{y}^{\mu}$ of $y$ under the action of $\mu (t)$ equals $y_{F(\mf{t})}$. 
Conversely, if $\mu (t) \subseteq T$ is a 1-PS such that $w=\hat{y}^{\mu}$, then there is a $c'\in \R$ such that $\mf{t}=\log (\mu)-c'{\bf 1}\in {\cal T}^+ (y)$ and $w=y_{F(\mf{t})}$.     
\end{enumerate} 
\end{lemma}

\noindent 
The proofs of all the above lemmas are straightforward. 

We now come to the main proposition of this subsection. We assume there is a 1-PS $\lambda(t) \subseteq T$, driving $y$ to $z$. Recall that for this $\lambda (t)$, we also have the 1-PS $\lambda '(t)$ (and an $\overline{\ell} \in {\cal T}$ with $t^{\overline{\ell}}=\lambda '(t) $) such that $z$ is the leading term of $y$ under $\lambda '$ of degree $0$.  This gives us an $\overline{\ell} \in {\cal T}^+ (y)$ such that $z=y_{F(\overline{\ell)}}$. 

\begin{prop} \label{prop:facedim}
Let $\lambda (t)\subseteq T$, $z=\hat{y}$ and $\overline{\ell}$ be as above. Suppose that $dim(\Xi (y))-dim(\Xi (z))\geq 2$, then there is a $\mf{t}\in {\cal T}^+ (y)$, its face $F=F(\mf{t})$ and a 1-PS $\mu (t)=t^{\mf{t}}$, such that (i) $y_F$ is the leading term of $y$ of degree $0$ under $\mu$ and $z$ is the leading term of $y_F$ under $\lambda '$ of degree $0$, and (ii) $dim(\Xi (y))> dim(\Xi (y_F )) > dim(\Xi (z))$.  Thus, there is an intermediate orbit $O(y_F )$ such that $dim(\Xi (y_F ))$ is strictly intermediate.
\end{prop}

\noindent 
{\bf Proof}: Let $F_y =\Xi (y)$ and suppose that $dim(F_y )\geq dim(F_z )+2$.  
Then, we must have $dim({\cal T}_z)\geq dim({\cal T}_y)+2$. 
By polyhedral theory, $\langle \chi , \overline{\ell} \rangle >0$ for all $\chi \in \overline{F}=F_y -F_z$. Thus $\overline{\ell} \in {\cal T}_z -{\cal T}_y$. 
Suppose that $\mf{s}$ is another element of ${\cal T}_z -{\cal T}_y$ which is linearly independent of ${\cal T}_y +\R \overline{\ell}$. Such an element exists since $dim({\cal T}_z )\geq dim({\cal T}_y )+2$. 
Let $\overline{a}=(a_{\chi})$ and $\overline{b}=(b_{\chi})$ be vectors defined on $\overline{F}$ such that for any $\chi \in \overline{F}$, $a_{\chi} =\langle \chi , \overline{\ell} \rangle$ and $b_{\chi}= \langle \chi, \mf{s} \rangle$. Clearly $\overline{a}>0$ and $\overline{b}$ is non-zero vector linearly independent of $\overline{a}$. 

Next, let us consider $\mf{t}(\epsilon)=\overline{\ell}+\epsilon \mf{s}$. Given the properties of $\overline{a}, \overline{b}$, there is an $\epsilon >0$ such that (i) $\langle \chi , \mf{t} (\epsilon) \rangle \geq 0$ for all $\chi \in \overline{F}$ {\em and} (ii) there is at least one $\chi ' \in \overline{F}$ for which $\langle \chi' , \mf{t}(\epsilon)\rangle =0$, and finally (iii) a $\chi'' \in \overline{F}$ such that $\langle \chi '' , \mf{t} (\epsilon) \rangle >0$.  We denote this by $\mf{t}'$. If $F=F_{\mf{t}'}$, then clearly $F_y \supsetneq F \supsetneq F_z$ and ${\cal T}_y \subsetneq {\cal T}_{y_F } \subsetneq {\cal T}_z$. 
Hence, for $\mu(t)=t^{\mf{t}'}$ we have $\hat{y}^{\mu}=y_F$ has the required property on dimensions that $dim(F_y ) > dim(F)>dim(F_z)$. Finally, it is easy to see that $\widehat{y_F}^{\lambda'}=z$. We come to the last part that $O(y_F)$ is intermediate to $O(y)$ and $O(z)$. Since $y_F$ is a limit of $y$ under $\mu$, we have $\overline{O(y)} \supseteq O(y_F )$, and since $z$ is a limit of $y_F$ under $\lambda$, we have $\overline{O(y_F)}\supseteq O(z)$. Thus $y_F$ is indeed an intermediate orbit. This proves the proposition. $\Box $

\begin{ex}
Let $X=\C^3, G=GL(X)$ and $V=Sym^2 (X^* )$. Let $B=\{ x_1 ,x_2 , x_3 \}$ be a basis for $X^*$. Let $y=x_1 x_2 +x_1 x_3 + x_2 x_3 + x_3^2$ and $\lambda(t)=diag(1,1,t)$. Thus $\ell =[0,0,1]$ and $t^{\ell}=\lambda (t)$. We have $z=x_1 x_2$ is the limit $\hat{y}^{\lambda}$. $H_0$ contains a torus $\mu(u)=u^{[1,-1,0]}=diag(u,u^{-1},1)$. Let $\ell ' =[1,-1,0]$. Note that the standard torus $T=diag(t_1 ,t_2 ,t_3)$ contains both $\lambda (t)$ and $\mu(u)$. We then choose it as the master maximal torus of $G$ and use it to construct $\Xi (V), \Xi (y)$ and $\Xi (z)$. We record $\Xi (y)=\{ [1,1,0], [1, 0, 1], [0,1,1], [0,0,2]\}$ and $\Xi (z)=\{ [1, 1, 0] \}$. We have ${\cal T}_y =0$ but ${\cal T}_z =\R \ell + \R \ell '$ has dimension $2$. Hence we may apply Prop. \ref{prop:facedim}. Indeed, we may choose $\ell '$ as the element $\mf{s}$ in the proof of the above proposition and construct $\mf{t} (\epsilon )= \ell + \epsilon \ell '$. We get $\epsilon=1$ and $\mf{t}'=[1,-1,1]$. We then see that $F=\{ [1,1,0], [0,1,1] \}, y_F =x_1 x_2 + x_2 x_3 $ and $O(y_F )$ is an intermediate limit.  
\end{ex}
\begin{corollary}
Let $z=\hat{y}$ under $\lambda (t)$ be such that there is no intermediate $G$-stable variety between $O(y)$ and $O(z)$. Then for any maximal torus $T$ containing $\lambda $, there is a $\lambda '' (t) \subseteq T,z''\in O(z)$ and $y'' \in O(y)$ such that $z''$ is the leading term of $y''$ under $\lambda '' (t)$ and $\Xi (z'')$ is a facet of $\Xi (y'')$. 
\end{corollary}

\noindent 
{\bf Proof}: We consider tuples $(y',z',\mu' )$, where $y' \in O(y), z'\in O(z), \mu'(t)\subseteq T$ and $z'=\hat{y'}^{\mu}$, and induct on $\delta (y',z',\mu )=dim(\Xi (y'))-dim(\Xi (z'))$. We begin with $(y,z,\lambda ')$, where $\lambda '(t)$ is the degree $0$ version of $\lambda(t)$, i.e., when $z$ is a limit of $y$ of degree $0$ under $\lambda '(t)$. If $dim(\Xi(y))=dim(\Xi(z))+1$, we are done. If not and $dim(\Xi (y))>dim(\Xi (z))+1$, then by the above proposition, there is a $\mu (t)\subseteq T$ and a $y_F$ whose orbit is intermediate and for which $dim(\Xi (y_F ))$ is strictly intermediate. Since, there is no intermediate variety between $y$ and $z$, we must either have (a) $\overline{O(y_F)}=\overline{O(y)}$ or (b) $\overline{O(y_F )}=\overline{O(z)}$. This implies that either (a) $y_F \in O(y)$ or (b) $y_F \in O(z)$. We thus get either (a) the tuple $(y_F , z, \lambda) $ or (b) the tuple $(y, y_F , \mu)$ (as the case may be), with a smaller $\delta$ and yet implementing the same limit. $\Box $

\begin{ex}
Let $X=\C^9$ be the space of $3\times 3$-matrices, and $G=GL(X)$. Let $V=Sym^3 (X^* )$ and $det_3 (X) \in V$ be $y$. The functions $B=\{ x_1 , \ldots ,x_9 \}$ is a basis for $X^*$ corresponding to the entries of the matrix. 
For any matrix $x\in X$, we have the expression $x=x_a+x_s$ decomposing $x$ as a sum of an antisymmetric and a symmetric matrices. For $\lambda_2 (t)(x)=tx_a +x_s$, we have $det_3 (\lambda_2 (t)x)=t Q_2(X) +t^3 R_2 (X)$. For a suitable choice of basis, we have:
\[ \lambda_2 (t) det_3 (X)=det\left( \left[
\begin{array}{ccc}
2tx_6 & tx_8 +x_1 & tx_9 -x_2 \\
tx_8 - x_1 & 2tx_5 & tx_7 + x_3 \\
tx_9 + x_2 & tx_7 - x_3 & 2tx_4 \\ \end{array} \right] \right) 
=tQ_2 (X) + t^3 R_2 (X)\]
where:
\[ \begin{array}{rcl}
Q_2 (X) &=&2( x_4 x_1^2 +x_5 x_2^2 +x_6 x_3^2 + x_7 x_1 x_2 +x_8 x_2 x_3 +x_9 x_1 x_3 )\\
R_2 (X) &=&8x_4 x_5 x_6 -2x_6 x_7^2 -2 x_4 x_8^2 -2 x_5 x_9^2 +2 x_7 x_8 x_9  
\end{array}
\]
We know that there is no intermediate orbit between $O(det_3 (X))$ and $O(Q_2 )$. Using the torus $T$ with the basis $x_1 ,\ldots , x_9$, we have:
\[ \Xi (Q_2)=  \begin{array}{|ccccccccc|} \hline 
x_1 & x_2 & x_3 & x_4 & x_5 & x_6 & x_7 & x_8 & x_9 \\ \hline \hline 
2 & 0 & 0 & 1 & 0 & 0 & 0 & 0 & 0 \\
0 & 2 & 0 & 0 & 1 & 0 & 0 & 0 & 0 \\
0 & 0 & 2 & 0 & 0 & 1 & 0 & 0 & 0 \\
1 & 1 & 0 & 0 & 0 & 0 & 1 & 0 & 0 \\
0 & 1 & 1 & 0 & 0 & 0 & 0 & 1 & 0 \\
1 & 0 & 1 & 0 & 0 & 0 & 0 & 0 & 1 \\ \hline 
\end{array}
\]
We may similarly build $\Xi (R_2)$. We see that $dim(\Xi (Q_2))=6, dim(\Xi (R_2))=4$ while $dim(\Xi (det_3 (X) )=7$. 
Thus $\Xi (Q_2 )$ is indeed a facet, i.e., a face of co-dimension $1$ within $\Xi(det_3 (X))$.

\end{ex}
\subsection{The Co-limit space $Z(\lambda)$  and $O(y)$} \label{subsec:CoLimit}

\noindent 
Let us resume our standard assumption of action by a 1-PS:
\[ \lambda (t) y= t^d y_d +t^e y_e +\ldots + t^D y_D \]
with the limit $z=\hat{y}$ and $y_e$ as the tangent of approach.

Let us act $\lambda$ on $y'$, a conjugate of $y$ to get:
\[ \lambda(t) \cdot y' = t^a y'_a + \cdots + t^b y'_b \]
Next, for any $d'$, we define $Y_{d'}$ and $Y$ as below:
\[ Y_{d'} =\{ y'=gy| y'_a =0 \mbox{ for all $a<d'$} \} \: \: \mbox{ and } Y=\cup_{d'} Y_{d'} \]
Thus, $Y_{d'} $ consists of those elements $y' \in O(y)$ for which $deg(\hat{y'})\geq d'$. Note that the notation and definitions are all with respect to this fixed $\lambda$. 

Let $V_{d'}$ be the degree $d'$ subspace of $V$ and consider the projection $\pi_{d'} : V \rightarrow V_{d'}$. We define $Z_{d'}=\pi_{d'} (Y_{d'})$ and $Z$ as $ \cup_{d'} Z_{d'}$. Thus, it is the space of all $z'$ which are degree-$d'$ limits under $\lambda $ of some conjugate $y'$ of $y$. Note that $y\in Y_d $ and $z\in Z_d $. We call $Z_d$ as the {\bf space of co-limits} of $z$.

The importance of $Z_d$ comes from the following lemma:

\begin{lemma}
Let $O(Z_d)=\{ gz' | z' \in Z_d \mbox{ and } g\in G\} $ and $\overline{O(Z_d )}$ be its closure. Then $\overline{O(Z_d)}$ is an intermediate variety, i.e, $\overline{O(z)}\subseteq \overline{O(Z_d )} \subseteq \overline{O(y)}$. 
\end{lemma}

\noindent 
{\bf Proof}: Since $z\in Z_d$, it is clear that $\overline{O(z)}\subseteq \overline{O(Z_d )}$. For the other inclusion, we see that every element $z'\in Z_d$ is the leading term $\hat{y'}$ for some $y' \in Y_d \subseteq O(y)$. Hence $z' \in \overline{O(y')}$. But $O(y')=O(y)$ and is $G$-stable and hence $\overline{O(z')} \subseteq \overline{O(y)}$. This proves the second inclusion. $\Box $

  All points in $\overline{O(Z_d)}$ will have a representative (up to conjugation) of pure degree $d$. If 
$y$ does not have this property, e.g., when $y$ is stable, then $\overline{O(Z_d )} \subsetneq \overline{O(y)}$. 

The more interesting question is whether $O(Z_d )$ contains $z'$ which are {\em not} in $O(z)$. To answer this, we will compare $T_z Z_d$ and $(T_z O(z))_d$, the degree $d$ component of the tangent space of $T_z O(z)$. Note that $(T_z O(z))_d \subseteq T_z Z_d $.  

Recall that we have the parabolic group $P(\lambda )$, its unipotent radical $U(\lambda)$ and a special reductive complement $L(\lambda)$. Recall also that the Lie algebra of $P(\lambda)$ is ${\cal P}(\lambda )=\sum_{i\geq 0} {\cal G}_i$, i.e., the subspace generated by elements within ${\cal G}$ of non-negative degree, and that of $L(\lambda)$ is ${\cal L}(\lambda )={\cal G}_0$. 
We have the following important lemma:

\begin{lemma} \label{lemma:YZequiv}
The map $\pi_d :Y_d \rightarrow Z_d$ is $P(\lambda )$ equivariant. 
\end{lemma}
\noindent
The proof is straightforward. The lemma implies that $L(\lambda )$ has an action on $Z_d$. 

\begin{assume}
We assume that both $Y_d$ and $Z_d$ are smooth at the points $y$ and $z$. See also Remark~\ref{assumption-smooth}   \end{assume}

\begin{lemma}Let $y \in Y_d$ and $z =\pi_d (y)\in Z_d$ be as above. Then we have:
\begin{enumerate} 
\item The tangent space $T_{z} Z_d$ is given by the image of $T_{y} Y_d$. In other words, $(\pi_d)_* T_{y} Y_d =T_{z} Z_d$. 
\item Let $O_{L(\lambda )} (z) \subseteq Z_d $ be the $L(\lambda)$-orbit of $z$. Let $T_{z} O(z)$ be the tangent space of the $G$-orbit of $z$ and $T_{z} O_{L(\lambda)} (z)$ be the tangent space for the $L(\lambda )$-orbit of $z$, then $T_{z} O_{L(\lambda)} (z)=(T_{z} O(z))_d ={\cal G}_0 z$.  
\end{enumerate} 
\end{lemma}

\noindent 
{\bf Proof}: The first part follows from the definition of $Z_d$ as the image of $Y_d $ the projection $\pi_d : V \rightarrow V_d $. For the second part, note that $T_{z} O(z) =\sum_i {\cal G}_i z$ and thus $(T_{z} O(z))_d ={\cal G}_0 z= {\cal L}(\lambda) z$. But this is precisely $T_{z} O_{L(\lambda)} (z)$. $\Box $

We illustrate the above lemma with an analysis of an example communicated to us by Professor V. Popov.

\begin{ex} \label{ex:popov}
Let $G=GL_4(\C)$ act by left multiplication on $V=\C^4 \oplus \C^4 \oplus \C^4$, represented as a $4\times 3$-matrix.
Thus, an element $v$ of $V$ may be viewed as a $4\times 3$ matrix and the action of a $g\in G$ is  
$g.v$ under the usual matrix multiplication. We will use $e_1, e_2, e_3$ and $e_4$ to denote the standard basis of $\C^4$ as column vectors.

We set $y = [e_1, e_2, e_3] \in V$ and $\lambda(t)=\mbox{diagonal}(1, t, t^2, t^2)$. Clearly,
$$\lambda(t). y = [e_1, 0, 0] + t [0, e_2, 0] + t^2 [0, 0, e_3] 
                = y_0 + t y_1 + t^2 y_2 $$
 Thus, $z=y_0=[e_1 ,0,0]$ is the leading term of $\lambda(t).y$ and the tangent of approach is $y_e=y_1 = [0, e_2, 0]$.
 The grading induced by $\lambda(t)$  on ${\cal G}$ as well as the weight-subspaces of $V$ of $\lambda(t)$ action are depicted in the following diagram.
 $$
\begin{array}{cc}
\left[ \begin{array}{c|c|cc}
0& -1 & -2 & -2 \\ \hline
1 & 0 & -1 & -1 \\ \hline
2 & 1& 0 & 0 \\
2 & 1 & 0 & 0\end{array} \right] 
& \left[ \begin{array}{ccc}
0& 0 & 0  \\ \hline
1 & 1 & 1 \\ \hline
2 & 2 & 2 \\
2 & 2 & 2\end{array} \right] \\
& \\
\mbox{Weight space for $gl_4$} & \mbox{Weight space for $V$} \end{array}
$$

Note that $O(y)$ is the collection of all matrices $y'$ of rank $3$ and $\overline{O(y)}$ equals $V$. Since $d=0$, we have $Y_0$ is the collection of all matrices $y'$  which have a non-zero first row. Thus $Z_0$ is the collection of all matrices $z'$ which are non-zero only in the first row. These give rise to an infinite collection of orbits consisting of the space of all rank $1$ matrices of which $z$ is one element. 
All these $z'\in Z_0$ are closely related to the point $z$, indeed, they have identical stabilizers, viz.,  $G_z =G_{z'}$ and their orbits have the same dimension, viz., $4$. 

On the other hand, it is easily checked that $W_{z,y_e}$ consists of all matrices where the third column is zero. Thus the $G$-stable spaces $O(Z_0)$ and $W_{z,y_e}$ are not comparable, and yet $\overline{O(z)}\subsetneq O(Z_0), W_{z,y_e} \subsetneq \overline{O(y)}$. 

Finally, note that $Y_1$ are all rank $3$ matrices with zero first row. $Z_1$ are all matrices which are non-zero only in the second row. But note that $O(Z_1)=O(Z_0 )$. 
\end{ex}

\noindent 
Let us now analyse $T_y Y_d $ and $T_z Z_d $. We begin with a definition:

\begin{defn}
Let $y, \lambda$ and $d$ be fixed as above. A element $\mf{g}\in {\cal G}$ is called a $d$-stabilizer iff $\mf{g}=\sum_{i} \mf{g}_i$ is such that $(\mf{g}y)_a =0$ for all $a<d$. Let ${\cal G}_{y,d}\subseteq {\cal G}$ be the collection of $d$-stabilizers of $y$. 
\end{defn}

In other words, ${\cal G}_{y,d}$ is the collection of ``Lie elements'' $g\in G$ such that $gy\subseteq Y_d$. 
The importance of ${\cal G}_{y,d}$ comes from the following lemma. 

\begin{lemma}
The tangent space $T_y Y_d$ is given by the set $ \{ \mf{g} y | \mf{g}\in {\cal G}_{y,d}\}$. 
\end{lemma}

\noindent 
{\bf Proof}: All elements in a neighborhood of $y\in Y_d$ are given by $\rho (e^{\mf{g}t})(y)$ for some $\mf{g}\in {\cal G}$. The tangent vector of the path $\beta_{\mf{g}}(t)=\rho (e^{\mf{g} t}) (y)$ at $y$ is precisely $\rho_1 (\mf{g})(y)$. That the path lies is $\oplus_{i\geq d} V_i$ is equivalent to the condition that $\mf{g} \in {\cal G}_{y,d}$. $\Box $

Let us compute $T_z Z_d $ and compare it with $(T_z O(z))_d ={\cal G}_0 z$. Since $T_z Z_d =(\pi_d)_* (T_y Y_d )$, we have $T_z Z_d =({\cal G}_{y,d}y)_d$, the degree $d$ component of the tangent space $T_y Y_d$. If $\mf{g}= \mf{g}_k+\ldots +\mf{g}_l \in {\cal G}_{y,d}$, then the action of $\mf{g}$ on $y$ gives us:

\[ \begin{array}{rcl} \mf{g}\cdot y &=& (\mf{g}_k +\ldots + \mf{g}_{k+(e-d)-1}+ \mf{g}_{k+(e-d)}+\ldots )(y_d +y_e + \ldots) \\
&=& \mf{g}_k y_d +\ldots \mf{g}_{k+(e-d)-1} y_d + (\mf{g}_k y_e + \mf{g}_{k+(e-d)} y_d )+\ldots+\ldots  \\
\end{array} \]
If $k>0$, then the leading term of $\mf{g}y$ has a degree greater than $d$, and does not contribute to $T_z Z_d$ and therefore is of no interest. If $k=0$, then the leading term is $\mf{g}_0 y_d \subseteq {\cal G}_0 z$ which is in $(T_z O(z))_d$, and thus does not cross the orbit of $z$. 
Thus the interesting situation is when $k<0$. Then the first $|k|$ terms of the above expression must vanish for it to contribute to $T_z Z_d$. It is the $(|k|+1)$-th term which gives us an element of $T_z Z_d$ which depends not only on $y_d$ but also on higher degree components of $y$. The analysis of these cancellations motivates us to study the structure of ${\cal G}_{y,d}$.  

Let us begin with three special nested families within ${\cal G}_{y,d}$ as specified below:

\begin{defn}
Let $\mf{g}=\sum_{i\geq a} \mf{g}_a$ be an element $\mf{g}_a $ of ${\cal G}_{y,d}$. We define $\hat{\mf{g}}$ as the leading term of $\mf{g}$. Let:
\[ {\cal G}^{\cal H}_{y,d} =\{  \hat{\mf{g}} |\mf{g}\in {\cal G}_{y,d} \mbox{ such that } \hat{\mf{g}}\in {\cal H} \} \]
Similarly, we define ${\cal G}_{y,d}^{{\cal H}_{\overline{y_e}}}$ (and resp. ${\cal G}_{y,d}^{\hat{\cal K}}$) as those $\mf{g}\in {\cal G}_{y,d}$ for which $\hat{\mf{g}}$ belongs to ${\cal H}_{\overline{y_e }}$ (resp. $\hat{\cal K}$). 
\end{defn}

Note that since the members of ${\cal G}^{\cal H}_{y,d}$ and other spaces are only the leading terms of elements of ${\cal G}_{y,d}$, these are graded subspaces of ${\cal G}$. Since ${\cal H}\supseteq {\cal H}_{\overline{y_e}} \supseteq \hat{\cal K}$, we have the corresponding containments ${\cal G}_{y,d}^{\cal H}\supseteq {\cal G}_{y,d}^{{\cal H}_{\overline{y_e}}} \supseteq {\cal G}_{y,d}^{\hat{\cal K}}$ and $dim({\cal G}_{y,d}^{\cal H}/{\cal G}_{y,d}^{{\cal H}_{\overline{y_e}}}) \leq dim ({\cal H}/{\cal H}_{\overline{y_e}})$ and $dim({\cal G}_{y,d}^{{\cal H}_{\overline{y_e}}}/{\cal G}_{y,d}^{\hat{\cal K}}) \leq dim ({\cal H}_{\overline{y_e}}/\hat{\cal K})$.

\begin{prop} \label{prop:lessthan}
Let $y,z$and $d$ be as above. 
Suppose that $dim({\cal H}/{\cal H}_{\overline{y_e}})=r$ and $dim({\cal H}_{\overline{y_e}}/\hat{\cal K})=s$. Then there is a subspace $F\subseteq {\cal G}_{y,d}$ of dimension at most $r+s$ and a basis $B=\{ \mf{g}_1 ,\ldots , \mf{g}_b \} \subset {\cal G}_{y,d}$ such that for all 
$\mf{g} \in {\cal G}_{y,d}$, $\mf{g} \in {\cal P}(\lambda)+{\cal K}+F$. 
Moreover the set $\hat{B}=\{ \hat{\mf{g}}_i | \mf{g}_i \in B\}$ is linearly independent in ${\cal G}_{y,d}^{\cal H}/ {\cal G}_{y,d}^{\hat{\cal K}}$.
\end{prop}

\noindent 
{\bf Proof}:  Let $B'=\{ \mf{g}_1 ,\ldots ,\mf{g}_b \}$ be elements of ${\cal G}_{y,d}$ such that $\hat{B'}= \{ \widehat{\mf{g}_i}|i=1,\ldots, b\} $ is a basis for ${\cal G}_{y,d}^{\cal H}/ {\cal G}_{y,d}^{\hat{\cal K}}$. Let $\mf{g}=\sum_{i\geq k} \mf{g}_i \in {\cal G}_{y,d}$ be an arbitrary element. We prove the assertion by induction on the degree of $\hat{\mf{g}}$, and on the containment of $\hat{\mf{g}}$ in the chain ${\cal G}\supseteq {\cal H}\supseteq {\cal H}_{\overline{y_e}} \supseteq \hat{\cal K}$.

Suppose that $\mf{g}$ is as above, i.e., $(\mf{g}y)_a=0$ for all $a<d$. We write the action of $\mf{g}$ on $y$ as before:
\[ \begin{array}{rcl} \mf{g}\cdot y &=& (\mf{g}_k +\ldots + \mf{g}_{k+(e-d)-1}+ \mf{g}_{k+(e-d)}+\ldots )(y_d +y_e + \ldots) \\
&=& \mf{g}_k y_d +\ldots \mf{g}_{k+(e-d)-1} y_d + (\mf{g}_k y_e + \mf{g}_{k+(e-d)} y_d )+\ldots+\ldots  \\
\end{array} \]

If $k\geq 0$ then $\mf{g} \in {\cal P}(\lambda)$ and the assertion holds. That takes care of the base case. Next consider the case when $k < 0$. We then have $\mf{g}_k y_d =0$. Thus $\hat{\mf{g}} \in {\cal H}$ and $\mf{g}\in ({\cal G}^{\cal H}_{y,d})_k$, the $d$-stabilizers with leading terms in ${\cal H}$ and of leading degree $k$. Then, by the choice of $B'$, there is an element $\mf{g}' \in \C \cdot F$ such that $\mf{g}-\mf{g}' \in {\cal G}^{\hat{\cal K}}_{y,d}$, i.e., where the leading term has dropped to $\hat{\cal K}$. 

We are then reduced to the case
where $\hat{\mf{g}}\in \hat{\cal K}$. In which case, there is an element $\mf{k}\in {\cal K}$ such that $\hat{\mf{k}}=\hat{\mf{g}}$. Since $\mf{k}\cdot y=0$, we have the element $\mf{g}'=\mf{g}-\mf{k}$ of ${\cal G}_{y,d}$ which has a leading term of degree greater than $k$. The set $B$ is constructed from $B'$ from the excess of ${\cal G}_{y,d}$ over ${\cal P}(\lambda)+{\cal K}$. This proves the proposition. 
$\Box$

\begin{remark}
If $y$ is in the null cone of $V$ for the $G$-action and $\lambda$ is the ``optimal'' 1-PS then ${\cal G}_{y,d}={\cal P}(\lambda)$, see \cite{hesselink1979desingularizations}, Lemma 4.6. Thus $F$ measures the deviation of $\lambda $ from the optimal 1-PS which drives $y$ to $0$.   
\end{remark}

\begin{prop}
\label{newproposed}
Let $TW_z = \pi_d (\{  \mf{g} \cdot y | \mf{g} \in F\})$ be the leading terms of $F\cdot y$, then $T_z Z_d$ =$TW_z + {\cal G}_0 z $. If $\overline{TW}_z$ denotes the quotient $(TW_z + {\cal G}_0 z)/{\cal G}_0 z$, then the codimension of $(T_z O(z))_d$ in $T_z Z_d $ is $dim(\overline{TW}_z)$.
\end{prop}

\noindent 
{\bf Proof}: Assuming manifold properties of $Y_d $ at $y$, we have \[\begin{array}{c}
dim(Y_d )=dim(T_y Y_d )\geq dim((T_y Y_d) _d )=dim(T_z Z_d ) =dim(Z_d )\\
dim((T_y Y_d )_d) = dim (({\cal G}_{y,d}y)_d) \end{array} \]
This gives us the equation $dim(({\cal G}_{y,d} y)_d)=dim(T_z Z_d)$. But, by Prop. \ref{prop:lessthan}, ${\cal G}_{y,d}={\cal P}(\lambda ) +{\cal K}+F$. Now $(({\cal P}(\lambda) +{\cal K})y)_d \subseteq (T_z O(z))_d$. 
Thus $F$ is the only part of ${\cal G}_{y,d}$ which leads to tangent vectors outside $(T_z O(z))_d$. Whence, $T_z Z_d = \pi_d (F)+ (T_z O(z))_d$. But $(T_z O(z))_d ={\cal G}_0 z$ hence $dim(T_z Z_d )-dim({\cal G}_0 z )$ is given by the expresion $dim(TW_z + {\cal G}_0 z ) -dim({\cal G}_0 z)$ and this is precisely $dim(\overline{TW}_z )$. This proves the proposition. 
 $\Box $

\begin{corollary}
If $dim(Z_d )>dim((T_z O(z))_d)$ then $F$ is non-empty. In other words, there is a $d$-stabilizer $\mf{g} \not \in {\cal P}(\lambda)+{\cal K}$. If $\lambda$
 were optimal then $\overline{O(z)}=\overline{O(Z_d )}$. 
\end{corollary}

\begin{remark}
\label{assumption-smooth}
The above computation of $T_z Z_d$ depends on $z\in Z_d$ being smooth. But is that so? If $z$ were smooth, then there would be the action of $H_0 =H \cap L(\lambda)$ on $T_z Z_d$ as follows. Since $\pi_d :Y_d \rightarrow Z_d $ is $L(\lambda)$-equivariant, and $H$ stabilizes $z$, we have the action of $H_0 $ on $Z_d $ which keeps $z$ fixed. This gives us an action of $H_0 $ and its Lie algebra, ${\cal H}_0$ on $T_z Z_d $. 

Even looking at the simple case when $\lambda $ has two components, $d=0$ and $e=1$, we have $T_z Z_d ={\cal H}_{-1} y_1 +{\cal G}_0 z $. But ${\cal G}_0 z$ is already preserved by $H_0$. Thus, we do have an $H_0$-action on $T_z Z_d$ if $H_0 ({\cal H}_{-1} y_1) \subseteq TW$. This seems to be a precondition for the smoothness of $z$ in $Z_d$.  

\end{remark}
\noindent 
We end this section with some examples.

\begin{ex} 
Continuing with Ex. \ref{ex:popov}, with $V=\C^{4 \times 3}, y=[e_1 , e_2 , e_3]$ and $\lambda(t)=\mbox{diagonal} (1, t, t^2, t^2)$:
$$\lambda(t). y = [e_1, 0, 0] + t [0, e_2, 0] + t^2 [0, 0, e_3] 
                = y_0 + t y_1 + t^2 y_2 $$
Thus $d=0, e=1, z=[e_1 ,0,0]$ and the tangent of approach is $y_1 =[0,e_2 ,0]$. The tangent space $T_z O(z)$ is $\{ [c,0,0]|  c \in \C^4 \}$, the space of all matrices which are $0$ in the second and third columns. Also recall the weight spaces:
$$
\begin{array}{cc}
\left[ \begin{array}{c|c|cc}
0& -1 & -2 & -2 \\ \hline
1 & 0 & -1 & -1 \\ \hline
2 & 1& 0 & 0 \\
2 & 1 & 0 & 0\end{array} \right] 
& \left[ \begin{array}{ccc}
0& 0 & 0  \\ \hline
1 & 1 & 1 \\ \hline
2 & 2 & 2 \\
2 & 2 & 2\end{array} \right] \\
& \\
\mbox{Weight space for $gl_4$} & \mbox{Weight space for $V$} \end{array}
$$
We see that: 
\[ {\cal K} = \hat{{\cal K}} =\left[ \begin{array}{cccc}0 & 0 & 0 & * \\
0 & 0 & 0 & * \\
0 & 0 & 0 & * \\
0 & 0 & 0 & *
\end{array} \right] 
\hspace*{0.4cm} {\cal H} =\left[ \begin{array}{cccc}
0 & * & * & * \\
0 & * & * & * \\
0 & * & * & * \\
0 & * & * & *
\end{array} \right] 
\hspace*{0.4cm} {\cal H}_{\overline{y_e}} =\left[ \begin{array}{cccc}
0 & 0 & * & * \\
0 & 0 & * & * \\
0 & 0 & * & * \\
0 & 0 & * & *
\end{array} \right] 
\]
Observe that ${\cal H} / {\cal H}_{\overline{y_e}}$ and ${\cal H}_{\overline{y_e}}/{\cal K}$ have
bases $\{ e_{1,2}, e_{2,2}, e_{3,2}, e_{4,2} \}$ and $\{e_{1,3}, e_{2,3}, e_{3,3}, e_{4,3}\}$ respectively.
It is easily verified that ${\cal G}_{y,d}$ equals ${\cal G}$, and is spanned by ${\cal P}(\lambda ) \cup {\cal K} \cup F$ where $F=\{ e_{1,2}, e_{1,3}, e_{2,3} \}$. Further,
$$ 
e_{1,2}.y = [0, e_1, 0], \;\; e_{1,3}.y = [0,0,e_1]\;\; e_{2,3}. y=[0,0,e_2]
$$
In the notation of Proposition~\ref{newproposed}, we see that $TW_0$, the vectors of weight $0$  are spanned by
$[0,e_1, 0]$ and $[0,0,e_1]$ and that $TW_0 =\overline{TW}_0$. These are elements of $T_z Z_0 - 
T_z O(z)$. Indeed, it is easily seen that $Z_0$ consists of all non-zero matrices where the only non-zero row is the first row. Consider for example, $z'=[e_1, \alpha_1 e_1, \alpha_2 e_1]$, where $\alpha_1 , \alpha_2 \in \C$. Then for $y'=[e_1 , \alpha_1 e_1 +e_2 , \alpha_2 e_1 +e_3 ]$, we see that $y' \in O(y)$ and $\widehat{y'}=z'$. Thus $z'\in Z_0$. 

Note that $O(z)_0 =\{[\alpha e_1, 0, 0]| \alpha \in \C, \alpha \neq 0\}$, and $O(z)$ are matrices of the form $[v,0,0]$ where
$v\neq 0$. However $O(Z_0)$ is the space of matrices of rank $1$. This has an infinite family of orbits, each of
the same dimension as $O(z)$ (i.e., $4$). Thus  $\overline{O(z)} \subsetneq \overline{O(Z_0 )} \subsetneq \overline{O(y)}$ is a strict intermediate variety.
\end{ex}

\begin{ex}
Consider the $G=GL_3 (\C)$-module $V=Sym^3 (\C^3)$ and the form $p(x,y,z)=x^2(x-y) + (x-2y)^2z+z^2(ax+by) +z^3 \in V$. We may choose $a,b$ so that ${\cal K}={\cal G}_p =0$. Let $\lambda (t)$ be such that $\lambda (t)x=x, \lambda (t)y=y$ but $\lambda (t)z=tz$. 
For the action of $\lambda (t)$ on $p$, we have $d=0, e=1$ and $\hat{p}=(x-y)x^2$ is the leading term and $q=(x-2y)^2 z$ is the tangent of approach $y_e$. Clearly ${\cal H}$, the stabilizer of $\hat{p}$ is $4$ dimensional and contains $Hom(\C z, \C \cdot \{ x,y,z\})$. The other basis element of ${\cal H}$ is the toric one parameter family of  $gl_3$ with $a_{11} = 1, a_{21}=3$ and $a_{22} =-2$, which is of degree $0$. 
However none of these elements stabilize $\overline{q}$ so ${\cal H}_{\overline{q}}=0$. Thus, $\hat{\cal K}={\cal H}_{\overline{q}} = 0$. The expected dimension of ${\cal W}_{z,y_e}$ is $dim(G) - dim({\cal H}_{\overline{y_e}}) = dim(G)$. Since this is also the dimension of $\overline{O(p)}$ there is no strictly intermediate variety of the type ${\cal W}_{z,y_e}$.

On the other hand the space $T_{\hat{p}} (O(\hat{p}))_0$ is given by ${\cal G}_0 \hat{p}$, which reduces to the action of $gl(x,y)$ on $\hat{p}$. This gives us the $4$ forms whose coefficients in terms of the basis elements in the first row of the matrix given below are:
\[ 
\begin{array}{|c||r|r|r|r|}\hline \hline 
& x^3 & x^2 y & xy^2 & y^3 \\ \hline \hline 
x\partial p/\partial x & 3 & -2 & 0 & 0 \\
y\partial p/\partial x & 0 & 3 & -2 & 0 \\ 
x \partial p/\partial y & -1 & 0 & 0 & 0 \\
y \partial p/\partial y & 0 & -1 & 0 & 0 \\ \hline 
\end{array}\]

The rank of the above matrix, and therefore the dimension of $(T_{\hat{p}} O(\hat{p}))_0$, is $3$.  
 $y\partial/\partial z$ is an element of ${\cal H}_{-1}=Hom(\C z , \C \{ x,y\})$ and is in ${\cal G}_{y,0}$. Applying this to $q$ gives us $y\partial q/\partial z =(x-2y)^2y \in ({\cal H}_{-1} q)_0 \in T_{\hat{p}} Z_0 $. But this element is not in $T_{\hat{p}} (O(\hat{p}))_0$. So $\overline{O(\hat{p}})\subsetneq \overline{O(Z_0 )}$. Now $Z_0$ contains all forms in $Sym^3 (\C \{ x,y\})$, the space of degree $3$ forms in the variables $x,y$. Thus, every form in $\overline{O(Z_0 )}$ is stabilized by a conjugate of $Hom(\C z , \C\{ x,y\})$. Therefore $\overline{O(Z_0 )} \subsetneq \overline{O(p)}$ which gives us a strict intermediate variety.

\end{ex}

\subsection{Alignment and Co-limits}

In this subsection, we connect the space $Z_0$ with the presence of alignments for the case when $V=Sym^n (X^*)$ and $ X=Y\oplus Z$ are the weight spaces for $\lambda $, as in Section \ref{subsec:weight-spaces}. Suppose that $f_0 \in Sym^n (Y^*)$ is obtained as a degree $0$ limit of some $f\in V$ under $\lambda (t)$. We show that there is alignment between $f$ and $f_0'$ for some co-limit $f_0' $ of $f_0$. 

As in Section \ref{subsec:weight-spaces}, let us suppose that there is a 1-PS $\lambda (t)\subseteq GL(X)$ such that $X=Y\oplus Z$, where $Y = X_0$ and $Z=X_1$, are the weight spaces. Suppose that:
\[ \lambda (t) f= t^0 f_0 + \ldots + t^n f_n \]
with $f_0 \in Sym^n (Y^*)$ as the leading term. As before, let $GL(X)_f =K$ and $  GL(X)_{f_0}=H$. Moreover, let $GL(Y)_{f_0}=H_Y$ be the restriction of $H$ to $Y$. 

Let ${\cal Z}=\{ z_1 ,\ldots ,z_s \}$ be a basis for $Z$ and ${\cal Y}=\{ y_1 ,\ldots ,y_r \}$ be a basis for $Y$ and $X=Y\oplus Z$. For simplicity, let us identify $Z$ with $\C^{s+r}$ treated as row vectors, and $\{ e_1 ,\ldots , e_s , e_{s+1}, \ldots e_{s+r}\}$ as the standard basis, with $e_i =z_i$, for $i=1,\ldots s$ and $e_{s+j}=y_j$ for $j>s$. 

Then every element $\mf{h}\in {\cal H}$, and $h\in H$ is of the form:
\[ \mf{h}=\left[ \begin{array}{cc}
\mf{a} & \mf{b} \\
0 & \mf{d} \end{array} \right] 
h=\left[ \begin{array}{cc}
a & b \\
0 & d \end{array} \right] 
\]
where $\mf{a}\in End(Z), \mf{b}\in Hom(Z,Y)$ and $\mf{d}\in {\cal H}_Y$ (and $a\in GL(Z), b\in End(Z,Y)$ and $d\in H_Y $). Recall also the subgroups $L(\lambda)$ and $P(\lambda)$. We will also use the unipotent subgroup $\overline{U}(\lambda )$, as shown below:
\[ \overline{U}(\lambda) =\left\{ \left[ \begin{array}{cc}
I_s & b \\
0 & I_r \end{array} \right] | b\in End(Z,Y) \right\}  
\]
Note that $\overline{U}(\lambda )\subseteq H$. 
 
We begin with a lemma:
\begin{lemma} \label{lemma:linalign}
Let $g\in GL(X)$ be a diagonalizable endomorphism. Then there is an $h\in \overline{U}(\lambda)$ such that $hgh^{-1} \in P(\lambda)$.     
\end{lemma}

Let us assume this for the moment. Then we have:
\begin{prop}
Let $g \in K$ be a semisimple element, then there is $u\in \overline{U}(\lambda )$ and a $u'\in U(\lambda)$ such that (i) $(g^u)^{u'}=g^{u'u}$ is an alignment between $f^{u'u}$ and $f_0'=\widehat{f^{u'u}}^{\lambda}$, and (ii)  $\widehat{f^{u'u}}^{\lambda}=\widehat{f^{u}}^{\lambda}$. Moreover, (iii) $f_0 ' \in Z_0$, is a co-limit of $f_0$ such that there is a common irreducible component $Z^i$ of $\overline{Z_0}$ which contains both $f_0$ and $f_0 '$. 
\end{prop}

\noindent 
{\bf Proof}: Let $g\in K$ be a semisimple element. By lemma $\ref{lemma:linalign}$, there is an $u\in \overline{U(\lambda)} \subseteq H$ such that $ugu^{-1}\in P(\lambda)$. Thus, $f^u$ satisfies the hypothesis of Prop. \ref{prop:unistablegeneral}, and thus, there is a $u'\in U(\lambda )$ such that $\widehat{f^u}$ and $(f^u)^{u'}=f^{u'u}$ have an alignment, viz., $g^{u'u}$, without changing the leading term $\widehat{f^u}^{\lambda}=f_0 '$. 

Coming to (iii), since $u\in \overline{U(\lambda)}$ is unipotent, there is an $X\in {\cal H}_{-1}$, such that $u=e^X$. Consider the 1-parameter {\em algebraic} family $f_0 (t)=\widehat{e^{tX}f} \subseteq Z_0$. This is an algebraic path connecting $f_0$ with $f_0 '=f_0(1)$. This implies that there must be a component $Z^i$ of $\overline{Z_0}$ containing both. This proves the proposition. $\Box$

\begin{remark} 
We see that (i) $f_0 (t)$ are co-limits of $f_0$ for all $t$, (ii) $f_0  '= f_0 (1)$ is aligned with $f^g$, a conjugate of $f$, and finally (iii) the derivative $f'(0)$ equals $Xf_1 \in T_{f_0} Z_0$
is as identified by Prop. \ref{newproposed}. 
The crux, of course is that $f_0 (t)$ for $t>0$ need not lie in the same orbit, so it is not clear that $f_0 \in \overline{O(f_0 ') }$. This leads us to the following conjecture.
\end{remark}

\begin{conjecture} \label{conj:stabilityalign}
If $f_0$ is $L(\lambda )$-stable and is the only form in $Sym^n (Y)$ with stabilizer $H_Y$, then there is an alignment between $f_0$ and a conjugate $f^g$ of $f$.    
\end{conjecture}

\noindent 

Let us come to the the proof of Lemma \ref{lemma:linalign}. Since $dim(Z)=s$ and $dim(Y)=r$, $P(\lambda )$ and $\overline{U(\lambda)}$ are of the form: 
\[ P(\lambda)=\left[ \begin{array}{cc} 
    A & 0 \\ 
    C & D \end{array} \right] \: \: \: \: 
    \overline{U(\lambda )}=
    \left[ \begin{array}{cc} 
    I_s & X \\ 
    0 & I_r \end{array} \right]
    \]
The proof then boils down to the linear algebra computation below:

\begin{lemma}
    Let $R\in \C^{(s+r)\times (s+r)}$ be a diagnolizable matrix, then there is a matrix $S$, of the form shown below, such that $W=SRS^{-1}$ is block lower triangular, i.e., of the form shown below:
    \[ S=\left[ \begin{array}{cc} 
    I_s & X \\ 
    0 & I_r \end{array} \right] \: \: \: \: 
    SRS^{-1}=W=
    \left[ \begin{array}{cc} 
    W_{11} & 0 \\ 
    W_{21} & W_{22} \end{array} \right]
    \]
\end{lemma}
\noindent
{\bf Proof}: Let $v_1 ,\ldots ,v_{r+s}$ be left eigenvectors of $R$, with $v_i R=\lambda_i v_i$ for some $\lambda_i \in \C$. Now, since $\{ v_1 ,\ldots ,v_{r+s}\}$ form a basis of $\C^{r+s}$, there are some $v_{i_1}\ldots, v_{i_s}$ such that $\{ v_{i_1},\ldots , v_{i_s}, y_i ,\ldots ,y_r\}$ is also a basis of $\C^{r+s}$. Let us assume that $i_1 =1, \ldots , {i_s}=s$. Let 
$B=[v_1 ,\ldots , v_s ,y_1 ,\ldots , y_r ]^T$. Then $B$ and $BRB^{-1}$ are of the form shown below:
\[ 
B=
    \left[ \begin{array}{cc} 
    A & B \\ 
    0& I_r \end{array} \right] \: \: \: \: 
BRB^{-1}=
    \left[ \begin{array}{cc} 
    D & 0 \\ 
    A'& B' \end{array} \right]
\]
where $D$ is the diagonal matrix $diag(\lambda_1 ,\ldots , \lambda_r)$. Let $C$ be as shown below. Then, by suitable operations, shown below we get:
\[ 
C=
    \left[ \begin{array}{cc} 
    A^{-1} & 0 \\ 
    0& I_r \end{array} \right] \: \: \: \: 
(CB)R(CB)^{-1}=
    \left[ \begin{array}{cc} 
    E' & 0 \\ 
    A''& B'' \end{array} \right] \: \: \: \: 
CB=
    \left[ \begin{array}{cc} 
    I_s & A^{-1}B \\ 
    0& I_r \end{array} \right]
\]
Thus, the required matrix is $CB$ above. $\Box$

\section{In Conclusion}

Geometric Complexity Theory (GCT) as an approach to key lower bounds problems in computational complexity theory, was proposed in \cite{mulmuley2001geometric}.  It is the study of specific forms (such as the determinant and the permanent) with distinctive stabilizers, their homogenized versions in different spaces and their orbit closures. It is also a belief that specific structures associated with their stabilizers will eventually lead us to the construction of obstructions which yield lower bounds. Our paper was an attempt to build a bridge between two specific approaches - the representation theoretic approach of \cite{mulmuley2001geometric, mulmuley2008geometric, burgisser2011overview} and others, and the geometric approach of \cite{mignon2004quadratic}, \cite{liemethodsGCT} and others. 

We have used the existence of a 1-PS $\lambda $ which must connect the implementation of the degree homogenized version of the permanent ($z$) as a determinant ($y$) as the starting point. This study has led us to stabilizer limits, alignment and weight spaces, normal cones, co-limits and other geometric structures which provide interesting insights into the problem. To the best of our knowledge the statement ${\hat{\cal K}} \subseteq {\cal H}_{\overline{y_e}} \subset {\cal H}$ is the first explicit connection made between the stabilizer $K$ of $y$ and the stabilizer $H$ of $z$ in $y$'s orbit closure. That allows us to probe the alignment of $\lambda $ with respect to the two stabilizers. In the case of the determinant versus permanent problem, the presence of alignment leads to combinatorial insights in both, the boundary of the orbit closure of the determinant and the ``implementation'' map $\phi$ of the permanent as a determinant. This has been explored in Section \ref{sec:leadingterm}. 

Indeed, the chain of Lie algebras $\hat{\cal K}\subseteq {\cal H}_{\overline{y_e}} \subset {\cal H}$ leads us to consider intermediate orbits between $\overline{O(z)}$ and $\overline{O(y)}$. The varieties  $Spec(\overline{J}_{z,y_e})$ (in the normal cone) and $Z_d$ are steps in that direction. If the intermediate subvariety conjecture is correct, it will allow us to construct a tower $(W_i)$ of intermediate varieties where each step is ``tight'', i.e., with equality between $\hat{\cal K}_i$ and $({\cal H}_i)_{\overline{y_{e_i}}}$. On the other hand, the structure of $Z_d$ appears to be connected with the presence of alignment. 

In summary, our approach allows a rich computational framework which connects the GCT approach to classical questions on orbit closures, local group action and stratification. 

There are several interesting questions which seem to lie on this path. We briefly list these. 

We know that the boundary of the orbit of the determinant is of codimension one. The techniques used here lead us to suspect that the components of this boundary arise from special 1-PS $\lambda$ which are well-aligned with $K$, its stabilizer. This is connected to finding matrix families preserved by a large subgroup of $K$. 

There are two stabilizer families which are of interest. The first is ${\cal G}_{y'}$ for $y' \in Y_d$, the stabilizers of elements within the orbit $O(y)$ which preserve the degree of the limit. The second, of course, is ${\cal G}_{\widehat{y'}}$, the stabilizers of the elements $z'=\widehat{y'} \in Z_d$.  Since $\widehat{{\cal G}_{y'}}\subseteq {\cal G}_{\widehat{y'}}$, this is also related to alignment in the vicinity of $z\in Z_d$. The invariance of $\hat{\cal K}$ under the action of $U(\lambda)$, puts a compact structure on the orbit space of stabilizer limits. The master result on 1-PS is of course that of Kempf \cite{kempf1978instability}. Hence the question of finding ``optimal'' 1-PS $\lambda$ which ``implement'' the limit $y \stackrel{\lambda}{\rightarrow} z$, if they exist, and their properties, seems important.  

Following Kempf, we may define the ``null cone'' of $z$ as the space ${\cal N}(z)=\{ v\in V |z\in \overline{O(v)} \}$. The stratification of ${\cal N}(z)$ too will yield important information on how stabilizers change. This connects with the construction of a ``local model'' of \cite{liemethodsGCT} which described the ${\cal G}$-action in a neighbourhood of $z$. 

The last two questions appear to be connected to the varieties $Spec(\overline{J_{z,y_e}})$ and $Z_d $ studied in this paper.

\bibliographystyle{alpha}
\bibliography{references}
\end{document}